%% file: JohYukTAC2012.tex
\documentclass[12pt,draftcls,onecolumn,fleqn]{IEEEtran}   

\usepackage{amsfonts}
\usepackage{amsmath}
\usepackage{amssymb}
\usepackage{caption}
\usepackage{epsfig}
\usepackage{graphicx}
\usepackage{latexsym}
\input{symbols}

\usepackage{enumerate}

\def\stp{{\cal T}}

\newlength{\noteWidth}
\long\def\notes#1{\ifinner{\tiny #1} \else
\marginpar{\parbox[t]{\noteWidth}{\raggedright\tiny #1}}
\fi\typeout{#1}}

\newtheorem{thm}{\bf{Theorem}}[section]

\newtheorem{lem}[thm]{\bf{Lemma}}

\newtheorem{defn}[thm]{\bf{Definition}}
\newtheorem{asn}[thm]{\bf{Assumption}}
\newtheorem{rem}[thm]{\bf{Remark}}

\begin{document}

\title{Stochastic Stabilization of Partially Observed and Multi-Sensor Systems Driven by Gaussian Noise under Fixed-Rate Information Constraints$^1$}

\author{Andrew P. Johnston$^2$ and Serdar Y\"uksel$^2$}
\maketitle
\footnotetext[1]{This paper is to appear in part at the IEEE Conference on Decision and Control, Hawaii, 2012.}
\footnotetext[2]{Department of Mathematics and Statistics, Queen's University, Kingston, Ontario, Canada, K7L 3N6. Research supported by the Natural Sciences and Engineering Research Council of Canada (NSERC). Email: a.johnston@queensu.ca, yuksel@mast.queensu.ca}

\begin{abstract}
We investigate the stabilization of unstable multidimensional partially observed single-sensor and multi-sensor linear systems driven by unbounded noise and controlled over discrete noiseless channels under fixed-rate information constraints. Stability is achieved under fixed-rate communication requirements that are asymptotically tight in the limit of large sampling periods. Through the use of similarity transforms, sampling and random-time drift conditions we obtain a coding and control policy leading to the existence of a unique invariant distribution and finite second moment for the sampled state. We use a vector stabilization scheme in which all modes of the linear system visit a compact set together infinitely often. We prove tight necessary and sufficient conditions for the general multi-sensor case under an assumption related to the Jordan form structure of such systems. In the absence of this assumption, we give sufficient conditions for stabilization.
\end{abstract}

\section{Introduction}

\subsection{Problem Statement} \label{sec:IntroPS}

In this paper, we consider the class of multi-sensor LTI discrete-time systems with both plant and observation noise. The system equations are given by
\begin{align} \mathbf{x}_{t+1} = \mathbf{A}\mathbf{x}_t + \mathbf{B}\mathbf{u}_t + \mathbf{w}_t, \quad \mathbf{y}^j_t = \mathbf{C}^j\mathbf{x}_t + \mathbf{v}^j_t, \quad 1 \leq j \leq M, \label{eq:multisystem} \end{align}
where $\mathbf{x}_t \in \mathbb{R}^n$ and $\mathbf{u}_t \in \mathbb{R}^m$ are the state and control action variables at time $t \in \mathbb{N}$ respectively. The observation made by sensor $j$ at time $t$ is denoted by $\mathbf{y}^j_t \in \mathbb{R}^{p_j}$. The matrices $\mathbf{A}$, $\mathbf{B}, \mathbf{C}^j$ and random vectors $\mathbf{w}_t, \mathbf{v}^j_t$ are of compatible size. The initial state, $\mathbf{x}_0$, is drawn from a Gaussian distribution.

\begin{asn} The noise processes $\{\mathbf{w}_t\}$ and $\{\mathbf{v}^j_t\}$ are each i.i.d. sequences of multivariate Gaussian random vectors with zero mean. At time $t$, both $\mathbf{w}_t$ and $\mathbf{v}^j_t$ are independent of $\mathbf{x}_t$ and each other. \end{asn}

\begin{figure}[h]
\centering
\includegraphics[scale=1]{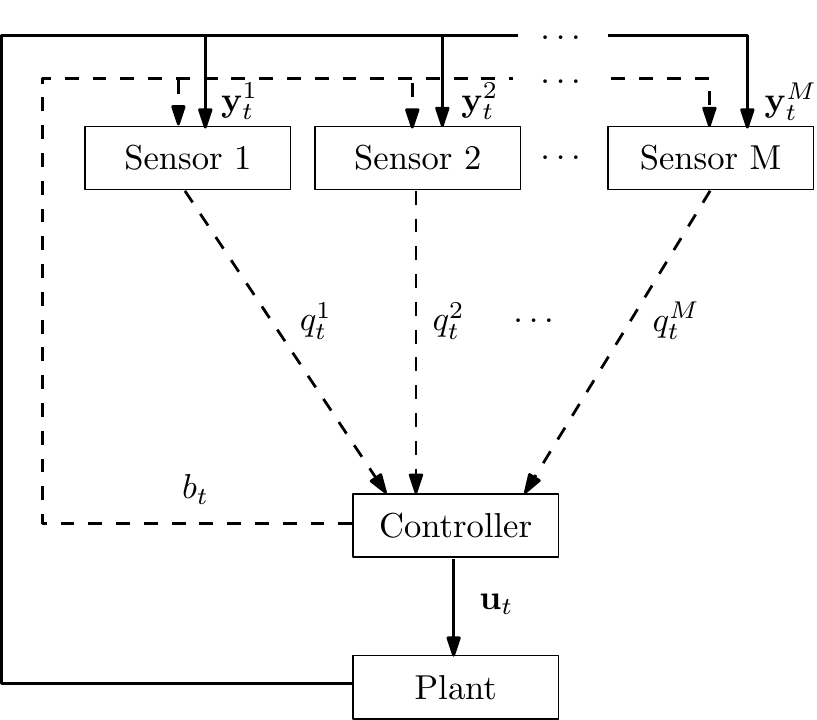}
\caption{A multi-sensor system with finite-rate communication channels.}
\label{fig:multisystem}
\end{figure}

\begin{asn} \label{ObservableAsn}
We require controllability and joint observability. That is, the pair $(\mathbf{A},\mathbf{B})$ is controllable and the pair $( [ (\mathbf{C}^1)^T \quad \cdots \quad (\mathbf{C}^M)^T ]^T, \mathbf{A} )$ is observable but the individual pairs $(\mathbf{C}^j, \mathbf{A})$ may not be observable.
\end{asn}
The setup is depicted in Figure \ref{fig:multisystem}. The observations are made by a set of $M$ sensors and each sensor sends information to the controller through a finite capacity channel. At each time stage $t$, we allow sensor $j \in \{1, \dots, M\}$ to send an encoded value $q^j_t \in \{1,2,\dots,N^j_t \}$ for some $N^j_t \in \mathbb{N}$. In addition, the controller can send a feedback value $b_t \in \{0,1\}$ at times $t=Ts$, where $T$ is the period of our coding policy and $s \in \mathbb{N}$. The value $b_t$ is seen by all sensors at time $t$. We define the rate at time $t$ as $ R_t = \sum^{M}_{j=1} \log_2(N^j_t)$. The coding scheme is applied periodically with period $T$ and so the rate for all time stages is specified by $\{N^j_0, \dots, N^j_{T-1} : 1 \leq j \leq M \}$. The average rate is
\begin{align} R_{\text{avg}} = \frac{1}{T} \left( M + \sum^{T-1}_{t=0} R_t \right), \label{eq:ravg} \end{align}
accounting for the encoded and feedback values.

\textbf{Information structure.} For a process $\{\mathbf{x}_t\}$ we define $ \mathbf{x}_{[a,b]} = \{\mathbf{x}_a, \mathbf{x}_{a+1},\dots, \mathbf{x}_b\}. $ At time $t$, each sensor $j$ maps its information $ I^{s_j}_t := \{ \mathbf{y}^j_{[0,t]}, b_{[0,t]} \} \to q^j_t \in \{1,\dots,N^j_t\}.$ The controller maps its information $ I^c_t := \{ q^1_{[0,t]}, \dots, q^M_{[0,t]} \} \to \mathbf{u}_t \in \mathbb{R}^m. $

\subsection{Notation} \label{sec:N}

We denote the indicator function of an event $E$ by $1_{E}$. We will use $\mathbb{R}^{m \times n}$ to denote the space of real $m \times n$ matrices and $\mathbb{R}^{n}$ to denote the space of real $n$ dimensional vectors. We let $\mathbb{R}^n_+$ be the space of real $n$ dimensional vectors with all entries nonnegative. Unless otherwise stated, all vectors are assumed to be column vectors. For any $\mathbf{x} \in \mathbb{R}^n$ we write $\mathbf{x} = \begin{bmatrix} x^1 & \cdots  & x^n \end{bmatrix}^T$ where $x^i \in \mathbb{R}$ is the $i^{th}$ entry. We define the absolute value operation for vectors as the component-wise absolute value. That is, $|\mathbf{x}| = \begin{bmatrix} |x^1| & \cdots  & |x^n| \end{bmatrix}^T$. For a matrix $\mathbf{A} \in \mathbb{R}^{n \times n}$ we denote its transpose by $\mathbf{A}^T$ and determinant by $\text{det}(\mathbf{A})$. If it is invertible, we denote the inverse by $\mathbf{A}^{-1}$. We let $\Lambda(\mathbf{A})$ denote the set of eigenvalues of $\mathbf{A}$. The $\ell^p$ norm is denoted by $\|\cdot \|_p$ and defined as $ \|\mathbf{x} \|_p = \left\{\sum^{n}_{i=1} |x^i|^p \right\}^{\frac{1}{p}}. $

\begin{defn} \label{vectorcomp} For $\mathbf{x} \in \mathbb{R}^n$ and $\mathbf{y} \in \mathbb{R}^n_+$ we write $\mathbf{x} \leq \mathbf{y}$ if $|x_i| \leq y_i$ for all $1 \leq i \leq n$. We write $\mathbf{x} \nleq \mathbf{y}$ otherwise.
\end{defn}

The observability matrix of sensor $j$ is $ \mathcal{O}_{(\mathbf{C}^j,\mathbf{A})} = \begin{bmatrix} (\mathbf{C}^j)^T & (\mathbf{C}^j \mathbf{A})^T & \cdots & (\mathbf{C}^j \mathbf{A}^{n-1})^T \end{bmatrix}^T, $ the null space is $ N^j = \text{Ker} (\mathcal{O}_{(\mathbf{C}^j,\mathbf{A})})$ and the observable subspace is defined to be $ O^j = (N^j)^{\perp}$ for $1 \leq j \leq M$.

\subsection{Brief Literature Review}

Due to space limitations, we are unable to give a fair account of the literature.  We refer the reader to the book \cite{MatveevSavkinBook} for a thorough review of the networked control literature and \cite{Yuksel2010} and \cite{Yuksel2010Tutorial} for a general overview of some of the related results.

There has been an extensive study in networked control theory regarding quantizer design for both stabilization and optimization. References \cite{WongBrockett}, \cite{Nair2004} and \cite{TatikondaMitter} obtained a lower bound on the average rate of the information transmission for the finiteness of second moments. For the system \eqref{eq:multisystem}, letting $\{\lambda_i\}$ be the set of eigenvalues of $\mathbf{A}$, this bound is $R_{\text{avg}} \geq R_{\min}$ where
\begin{align} R_{\min} = \sum_{|\lambda_i|>1} \log_2(|\lambda_i|). \label{eq:rmin} \end{align}

Various publications have studied the characterization of minimum information requirements for multi-sensor and multi-controller linear systems with an arbitrary topology of decentralization and the fundamental bounds have been extensively studied in \cite{MatveevSavkinBook}, \cite{EliaMitter} \cite{Brockett2000}, \cite{Tatikonda2003}, \cite{Nair04}, \cite{Nair07}, \cite{Savkin04}, \cite{Savkin05}, \cite{Gupta}, \cite{YukBasTAC072} and \cite{MatveevSavkin08}.

When a linear system is driven by unbounded noise, the analysis is particularly difficult since the bounded quantizer range leads to a transient state process (see Proposition 5.1 in \cite{Nair2004} and Theorem 4.2 in \cite{YukselBasarTAC2011}). For such a noisy setup, a stability result of the form  $\limsup_{t \to \infty} E[\|\mathbf{x}_t\|_2] < \infty$ was given for noisy systems with unbounded support in \cite{Nair2004}, which uses a variable-rate quantizer. Under this scheme, the quantizer is applied with a very high rate during some time intervals. More recently, a fixed-rate scheme was presented in \cite{Yuksel2010} for a scalar noisy system using martingale theory, which achieved the lower bound plus an additional symbol required for encoding. The existence of an invariant distribution was established under the coding and control policy presented, along with a finite second moment of the state. That is, $\lim_{t \to \infty} E[\|\mathbf{x}_t\|_2] < \infty$. \cite{YukMeynTAC2010} considered a general random-time stochastic drift criteria for Markov chains and applied it to binary erasure channels in a similar spirit.

\subsection{Contributions}
In view of the literature, the contributions of this work are as follows:
\begin{itemize}
\item The case where the system is multi-dimensional and driven by unbounded noise over a discrete-channel has not been studied to our knowledge, regarding the existence of an invariant distribution and ergodicity properties. Results for the limit properties of the finite moment are also new.
\item We give sufficient conditions for multi-sensor systems with both system noise and observation noise with unbounded support, which has not been treated previously, to our knowledge.
\end{itemize}
Our approach builds on the martingale and the random-drift programs considered in \cite{Yuksel2010} and \cite{YukMeynTAC2010}, however, new geometric constructions are needed for the vector and partially observed settings. We define a more general class of stopping times and adopt a further geometric approach.

We structure the paper as follows. In Section \ref{sec:SSS}, we study single-sensor systems and give our main result for such systems, Theorem \ref{MainResult}. Section \ref{sec:OPSSS} outlines the proof of Theorem \ref{MainResult}. The more detailed proofs can be found in Section \ref{sec:SRSSS}. In Section \ref{sec:MSS}, we study multi-sensor systems and give our main result for such systems, Theorem \ref{MultiSensorResult}. A supporting proof can be found in Section \ref{sec:PMSR}. Some basic definitions and results from the theory of matrix algebra, Markov chains and stochastic stabilization are provided in Section \ref{sec:MCSC}.

\section{Single-Sensor Systems} \label{sec:SSS}

\subsection{Problem Statement} \label{sec:PS}

Consider the class of single-sensor LTI discrete-time systems with both plant and observation noise. The system equations are given by
\begin{align} \mathbf{x}_{t+1} =  \mathbf{A} \mathbf{x}_t + \mathbf{B} \mathbf{u}_t + \mathbf{w}_t, \quad \mathbf{y}_t = \mathbf{C} \mathbf{x}_t + \mathbf{v}_t, \label{eq:system}  \end{align}
where $\mathbf{x}_t \in \mathbb{R}^n$, $\mathbf{u}_t \in \mathbb{R}^m$ and $\mathbf{y}_t \in \mathbb{R}^p$ are the state, control action and observation at time $t$ respectively. The matrices $\mathbf{A}, \mathbf{B}, \mathbf{C}$ and the noise vectors $\mathbf{w}_t, \mathbf{v}_t$ are of compatible size. The initial state, $\mathbf{x}_0$, is drawn from a Gaussian distribution. We label the eigenvalues of $\mathbf{A}$ as $\lambda_1,\dots,\lambda_n$.
Without loss, we assume that $\mathbf{A}$ is in real Jordan normal form and that $|\lambda_i| > 1$ for all $1 \leq i \leq n$.

\begin{asn} The noise processes $\{\mathbf{w}_t\}$ and $\{\mathbf{v}_t\}$ are each i.i.d. sequences of multivariate Gaussian random vectors with zero mean. At time $t$, both $\mathbf{w}_t$ and $\mathbf{v}_t$ are independent of $\mathbf{x}_t$ and eachother. \end{asn}

\begin{asn} The pair $(\mathbf{A},\mathbf{B})$ is controllable and the pair $(\mathbf{C},\mathbf{A})$ is observable.  \end{asn}

\begin{figure}[h]
\centering
\includegraphics[scale=1]{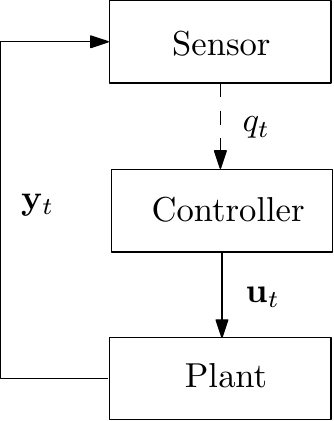}
\caption{A single-sensor system with finite-rate communication channel.}
\label{fig:system}
\end{figure}

The setup is depicted in Figure \ref{fig:system}. The observations are made by the sensor and sent to the controller through a finite capacity channel. At each time stage $t$, we allow the sensor to send an encoded value $q_t \in \{1,\dots,N_t\}$ for some $N_t \in \mathbb{N}$. We define the rate of our system at time $t$ as $ R_t = \log_2(N_t). $ Now, suppose that the channel is used periodically, every $T$ time stages. The rate for all time stages is then specified by $\{N_0, \dots, N_{T-1}\}$. The average rate is
\begin{align} R_{\text{avg}} = \frac{1}{T} \sum^{T-1}_{t=0} R_t. \label{eq:ssravg} \end{align}

\textbf{Information structure.}  At time $t$, the sensor maps its information $ I^s_t := \{ \mathbf{y}_{[0,t]} \} \to q_t \in \{1,\dots,N_t\}.$ The controller maps its information $ I^c_t := \{ q_{[0,t]}\} \to \mathbf{u}_t \in \mathbb{R}^m. $

\subsection{Main Result} \label{sec:SSMR}
Our main result for single-sensor systems is the following:
\begin{thm} \label{MainResult}
There exists a coding and control policy with average rate $R_{\text{avg}} \leq 1/(T2n)$ $\sum^{n}_{i=1} \log_2( \lceil |\lambda_i|^{T2n} + \epsilon \rceil + 1) $ for some $\epsilon > 0$ which gives:
 \begin{enumerate}[(a)]
\item the existence of a unique invariant distribution for $\{\mathbf{x}_{2nt}\}$;
\item $ \lim_{t \to \infty} E[\|\mathbf{x}_{2nt}\|_2] < \infty. $
    \end{enumerate}
\end{thm}

\begin{thm} \label{MinimumRate} The average rate in Theorem \ref{MainResult} achieves the minimum rate \eqref{eq:rmin} asymptotically for large sampling periods. That is, $\lim_{T \to \infty} R_{\text{avg}} = R_{\min}$. \end{thm}

\subsection{Coding and Control Policy} \label{sec:CCP}


For now, assume that $\mathbf{A}$ has only one eigenvalue $\lambda$. We will see later how this assumption can be made without loss.

Put $K = \lceil |\lambda| + \epsilon \rceil$ for some parameter $\epsilon > 0$ and consider the following scalar $(K+1)$-bin uniform quantizer. Assuming that $K$ is even, this is defined for $k \in \{1,2,\dots,K\}$ as
$$ Q^{\Delta}_K(x) = \begin{cases} \left( \frac{-(K+1)}{2} +k \right)\Delta, & \text{if $x \in [ \left( \frac{-K}{2} + k-1 \right)\Delta, \left( \frac{-K}{2} + k \right)\Delta )$},   \\
\frac{K-1}{2} \Delta, &  \text{if $|x| = \frac{K}{2} \Delta $},  \\ 0,  &  \text{if $|x| > \frac{K}{2} \Delta$},  \end{cases} $$
where $\Delta \in \mathbb{R}_+$ is the bin size. The set $[-\frac{K}{2}\Delta,\frac{K}{2}\Delta]$ is called the granular region while the set $(-\infty,-\frac{K}{2}\Delta) \cup (\frac{K}{2}\Delta, \infty)$ is called the overflow region. If the state is in the granular region, that is if $|x| \leq \frac{K}{2} \Delta$ then we say the quantizer is \textit{perfectly-zoomed}. Otherwise, we say it is \textit{under-zoomed}.

We write our quantizer as the composite function $Q^{\Delta}_K(x) = \mathcal{D}^{\Delta}_{K}(\mathcal{E}^{\Delta}_{K}(x))$. The encoder $\mathcal{E}^{\Delta}_{K} : \mathbb{R} \to \{1,\dots,K+1\}$ and decoder $\mathcal{D}^{\Delta}_{K} : \{1,\dots,K+1\} \to \mathcal{C}$ for $k \in \{1,2,\dots,K+1\}$ are
\begin{align} \mathcal{E}^{\Delta}_{K}(x) = \begin{cases} k, & \text{if  $x \in [ ( \frac{-K}{2} + k - 1 )\Delta,$} \\
 & \quad \quad ( \frac{-K}{2} + k )\Delta ), \\ K, & \text{if $x = \frac{K}{2} \Delta$}, \\ K +1, & \text{if $|x| > \frac{K}{2} \Delta$}, \\  \end{cases}
\mathcal{D}^{\Delta}_{K}(x) = \begin{cases} \left( -\frac{K+1}{2} + x \right)\Delta, \\ \qquad \qquad \text{if $x \neq K+1$},  \\ 0, \\ \qquad \qquad \text{if $x = K + 1$.}   \end{cases} \nonumber  \end{align}
At time $t$, we associate with each component $x^i_t$ a bin size $\Delta^i_t$. Let $q^i_t = \mathcal{E}^{\Delta^i_t}_{K}(y^i_t)$. We will be applying our control policy to system \eqref{eq:relabel} where $\mathbf{y}_s$ is a meaningful estimate of the state $\mathbf{x}_s$. Let our fixed rate be $N_t = K^n + 1$ for all $t \in \mathbb{N}$. Choose any invertible function $f: \{1,\dots,K\}^n \to \{1,\dots,K^n\}$. We then choose the encoded value
$$ q_t = \begin{cases}  f(q^1_t,\dots,q^n_t), & \text{if $q^i_t \neq 0$ for all $1 \leq i \leq n$}, \\
0, & \text{otherwise}. \end{cases}$$
Upon receiving $q_t \neq 0$, the controller knows $q^1_t,\dots, q^n_t$. The controller forms the estimate $\hat{\mathbf{x}}_t$ as $\hat{\mathbf{x}}_t = \begin{bmatrix} \hat{x}^1_t & \cdots & \hat{x}^n_t \end{bmatrix}^T, $ where
$$ \hat{x}^i_t = \begin{cases}\mathcal{D}^{\Delta^i_t}_{K}(q^i_t)  , & \text{if $q_t \neq 0$}, \\
0, & \text{otherwise}.   \end{cases} $$
We assume without loss that $\mathbf{A}$ is a Jordan block with eigenvalue $\lambda$. From the real Jordan canonical form (see for example \cite{Horn}), we know that it can be written as
$$ \mathbf{A} =  \begin{bmatrix} \lambda & 1 \\ & \lambda & \ddots \\ & &  \ddots & 1 \\ & & & \lambda \end{bmatrix}, \quad \text{if $\lambda \in \mathbb{R}$},  \quad
\mathbf{A} = \begin{bmatrix} \mathbf{D} & \mathbf{I} \\ & \mathbf{D} & \ddots \\ & &  \ddots & \mathbf{I} \\ & & & \mathbf{D} \end{bmatrix},  \quad \text{if $\lambda \in \mathbb{C}$},  $$
where in the complex case we write $\lambda = a + ib$ for some $a,b \in \mathbb{R}$ and define
$$ \mathbf{D} = \begin{bmatrix} a & b \\ -b & a \end{bmatrix}. $$
The update equations are
\begin{align} \Delta_{t+1} = \bar{Q}\left( q_t, \Delta_t \right) \Delta_{t}, \quad
\bar{Q}\left( q_t, \Delta_t \right) = \begin{cases} \rho |\lambda|, & \text{if $q_t = 0$}, \\
 \beta(\Delta_t), &\text{otherwise}, \end{cases} \label{eq:update} \end{align}
for some $\rho > 1$ and with
\begin{align} \beta(\Delta_t) = \text{diag}(\beta_1(\Delta^1_t),\dots, \beta_n(\Delta^n_t) ), \quad
 \beta_i(\Delta^i_t) = \begin{cases} 1, & \text{if $\Delta^i_t \leq L^i$ }, \\
 \frac{|\lambda|}{|\lambda| + \epsilon - \eta }, & \text{otherwise}, \end{cases} \label{eq:beta}
  \end{align}
for some $0 < \eta < \epsilon$ and $\mathbf{L} \in \mathbb{R}^n_+$. Note that if we define $ \bar{\mathbf{L}} = \mathbf{L} |\lambda| / (|\lambda| + \epsilon - \eta) $ then $ \Delta^i_t > \bar{L}^i $ for all $1 \leq i \leq n$ and all $t \in \mathbb{N}$.

\textbf{Bin ordering. } We set $ \mathbf{L} = c \Delta_0$, for some $0 < c \leq 1$. First let $\lambda \in \mathbb{R}$. For any $\delta > 0$ we can choose $\Delta^i_0$ and $\Delta^{i+1}_0$ such that $ \Delta^{i+1}_0 \leq \delta \Delta^i_0$ for all $1 \leq i \leq n-1$. With our update equations and our choice of $\mathbf{L}$ we get that the ordering is preserved over all time stages. That is, $ \Delta^{i+1}_t \leq \delta \Delta^i_t$ for all $1\leq i \leq n-1$ and $t \in \mathbb{N}$.

Now let $\lambda \in \mathbb{C}$. We choose $ \Delta^i_0 = \Delta^{i+1}_0$ for all $i$ odd. Thus, we have divided the complex modes into their conjugate pairs and set their initial bin sizes to be equal.  Our initial condition implies that $ \Delta^i_t = \Delta^{i+1}_t  $ for all $i$ odd and $t \in \mathbb{N}$. For any $\delta > 0$ we can choose $\Delta^i_0$ and $\Delta^{i+2}_0$ such that $\Delta^{i+2}_t \leq \delta \Delta^i_t$ for all  $1 \leq i \leq n-2$ and $t \in \mathbb{N}$.

Under our information structure, the update equations \eqref{eq:update} can be applied at the sensor and the controller. Our vector quantizer is implementable and at time $t$ the controller knows $\hat{\mathbf{x}}_t$. We choose the control action
$\mathbf{u}_t = -\mathbf{A} \hat{\mathbf{x}}_t.$

\subsection{Outline of Proof for Theorem \ref{MainResult}} \label{sec:OPSSS}

In this section, we outline the supporting results and key steps in proving our main result for single-sensor systems, Theorem \ref{MainResult}.

\begin{lem} \label{BeliefLemma} We can sample every $2n$ time stages and apply a similarity transform to $\mathbf{x}_t$ in \eqref{eq:system}  to obtain $\bar{\mathbf{x}}_s = \mathbf{P} \mathbf{x}_{2ns} $ with $s \in \mathbb{N}$ for some invertible matrix $\mathbf{P}$. This new state satisfies the following system of equations:
\begin{align} \bar{\mathbf{x}}_{s+1} = \bar{\mathbf{A}} \bar{\mathbf{x}}_s + \bar{\mathbf{u}}_s + \bar{\mathbf{w}}_s, \quad
\bar{\mathbf{y}}_s = \bar{\mathbf{x}}_s + \bar{\mathbf{v}}_s. \label{eq:sample} \end{align}
 The control action $\bar{\mathbf{u}}_s \in \mathbb{R}^n$ is chosen arbitrarily by the controller and the elimination of the $\mathbf{B}$ matrix can be justified by sampling. The estimate $\bar{\mathbf{y}}_s \in \mathbb{R}^n$ at time $s$ is known by the sensor. The noise processes $\{\bar{\mathbf{w}}_s\}$ and $\{\bar{\mathbf{v}}_s\}$ are each i.i.d. sequences of zero mean multivariate Gaussian random vectors. At time $s$, $\bar{\mathbf{w}}_s$ and $\bar{\mathbf{v}}_s$ are independent of $\bar{\mathbf{x}}_s$ but may be correlated with eachother. For $s_1 \neq s_2$, the vectors $\bar{\mathbf{w}}_{s_1}$ and $\bar{\mathbf{v}}_{s_2}$ are independent. The matrix $\bar{\mathbf{A}}$ is in real Jordan normal form and has eigenvalues $\lambda^{2n}_1,\dots,\lambda^{2n}_n$.
\end{lem}

By a slight abuse of notation, we will rewrite system \eqref{eq:sample} as
\begin{align} \mathbf{x}_{s+1} = \mathbf{A} \mathbf{x}_s + \mathbf{u}_s + \mathbf{w}_s, \quad
\mathbf{y}_s = \mathbf{x}_s + \mathbf{v}_s, \label{eq:relabel} \end{align}
where $\mathbf{x}_s \in \mathbb{R}^n$,  $\mathbf{u}_s \in \mathbb{R}^n$ and $\mathbf{y}_s \in \mathbb{R}^n$ are the state, control action and observation at time $s$ respectively. 

\begin{rem} We consider the case where $\mathbf{A}$ is a single Jordan block with eigenvalue $\lambda$. We can do this without loss since we are considering the single-sensor case and the sensor obtains an estimate for all components, as seen in Lemma \ref{BeliefLemma}. Thus, we can simply apply our control policy to each Jordan block. In all remaining theorems of this section, we will work with system \eqref{eq:relabel}. Where necessary, we will distinguish between the real and complex eigenvalue cases. 
\end{rem}

\begin{lem} \label{Markov} The process $\{(\mathbf{x}_s, \Delta_s)\}$ is Markov.
\end{lem}

Section \ref{sec:CCP} gives our control policy in terms of the parameters $\rho, \epsilon$ and $\eta$.

\begin{lem} \label{Irreducible}
For appropriate choices of $\rho, \epsilon$ and $\eta$, we can form a countable state space $\mathcal{S}$ for $\{\Delta_s\}$. The process $\{(\mathbf{x}_s,\Delta_s)\}$ is an irreducible Markov chain on $\mathbb{R}^n \times \mathcal{S}$.
\end{lem}

Define the sequence of stopping times
\begin{align*} \tau_0 = 0, \quad \tau_{z+1} = \min \left \{s > \tau_z : |\mathbf{y}_s| = |\mathbf{x}_s + \mathbf{v}_s| \leq \frac{K}{2}\Delta_s \right \}. \nonumber \end{align*}
These are the times when all quantizers are perfectly-zoomed. We assume that this is satisfied at time $s=0$. This technical condition is justified by showing that the process $\{(\mathbf{x}_s, \Delta_s)\}$ moves to such a perfectly zoomed state in a finite time, which is dominated by a geometric distribution (see a similar discussion in \cite{YukMeynTAC2010}).

\begin{thm} \label{GeometricDecay}
If $K$ is even then the following hold.
\begin{enumerate}[(a)]
\item For any $r > 0$ and any polynomial of finite degree $Q(k)$ there exists a sufficiently large $H$ such that $ Q(k) P(\tau_{z+1} - \tau_z > k \mid \mathbf{x}_{\tau_z},\Delta_{\tau_z}) \leq r^{-k}$ for all $k>H$ and for all $z \in \mathbb{N}$.
\item Let $\Delta_{\tau_z} \to \infty$ be equivalent to stating that $\Delta^i_{\tau_z} \to \infty$ for all $1 \leq i \leq n$. Then
$$ \lim_{\Delta_{\tau_z} \to \infty} P(\tau_{z+1} - \tau_z > 1 \mid \mathbf{x}_{\tau_z},\Delta_{\tau_z}) = 0 $$
uniformly in $\mathbf{x}_{\tau_z}$.
\end{enumerate}
\end{thm}

We define the compact sets
\begin{align*} S &= S_{\mathbf{x}} \times S_{\Delta}, \quad S_{\Delta} = \{ \Delta \in \mathbb{R}^n_+ : \Delta^i \leq F, \text{$1 \leq i \leq n$} \}, \\
& S_{\mathbf{x}} = \left\{\mathbf{x} \in \mathbb{R}^n : |x^i|  \leq \frac{K}{2} F, \text{$1 \leq i \leq n$}  \right\}, \end{align*}
for some $F > L^1$ where $L^1$ is a component of $\mathbf{L}$ as described in Section \ref{sec:CCP}. Note that at the stopping time $\tau_z$, if $\Delta_{\tau_z} \in S_{\Delta}$ then $ |x^i_{\tau_z}| \leq \frac{K}{2} \Delta^i_{\tau_z} \leq \frac{K}{2} F$, for all $1 \leq i \leq n$,
and thus $\mathbf{x}_{\tau_z} \in S_{\mathbf{x}}$ and $(\mathbf{x}_{\tau_z},\Delta_{\tau_z}) \in S$.

\begin{lem} \label{Drift}
For some $\gamma > 0$, the following drift condition holds:
\begin{align} & \gamma E \left[ \left. \sum^{\tau_{z+1} -1}_{s=\tau_z} (\Delta^1_s)^2 \right \vert \mathbf{x}_{\tau_z}, \Delta_{\tau_z}\right]  \leq (\Delta^1_{\tau_z})^2 - E[(\Delta^1_{\tau_{z+1}})^2 \mid \mathbf{x}_{\tau_z}, \Delta_{\tau_z}] + b1_{ \{ (\mathbf{x}_{\tau_z}, \Delta_{\tau_z}) \in S \}}. \label{eq:drift} \end{align}
For $\lambda \in \mathbb{C}$, the above also holds with $\Delta^2$ in place of $\Delta^1$.
\end{lem}



For $\mathbf{x} \in \mathbb{R}^n$, we say that $x^i$ and $x^{i+1}$ are a \textit{conjugate pair} if $i$ is odd. To simplify notation in the complex eigenvalue case we find it convenient to define for any $\mathbf{x} \in \mathbb{R}^n$, the set of vectors
$$\tilde{\mathbf{x}}^i = \begin{bmatrix} x^i & x^{i+1} \end{bmatrix}^T, \quad \text{if $i$ is odd}, \quad \tilde{\mathbf{x}}^i = \begin{bmatrix} x^{i-1} & x^{i} \end{bmatrix}^T, \quad \text{if $i$ is even}, $$
for $1 \leq i \leq n$. Note that $\tilde{\mathbf{x}}^i = \tilde{\mathbf{x}}^{i+1}$ for $i$ odd. We are only concerned with the case when $n$ is even.

\begin{thm} \label{MomentBound}
Let $\lambda \in \mathbb{R}$. For $i=n$, there exists a $\kappa > 0$ such that
\begin{align}  E\left[ \left. \sum^{\tau_{z+1} -1}_{s=\tau_z} (x^i_s)^2 \right \vert \mathbf{x}_{\tau_z}, \Delta_{\tau_z} \right] \leq \kappa (\Delta^1_{\tau_z})^2. \label{eq:MomentBound}  \end{align}
If $ \lim_{s \to \infty} E[(x^k_s)^2 ] < \infty $ then the above holds for $i = k-1$.

For $\lambda \in \mathbb{C}$, with $i=n-1$, there exists a $\kappa > 0$ such that
$$   E\left[ \left. \sum^{\tau_{z+1} -1}_{s=\tau_z} ( \tilde{\mathbf{x}}^i_s)^T \tilde{\mathbf{x}}^i_s \right \vert \mathbf{x}_{\tau_z}, \Delta_{\tau_z} \right] \leq \kappa (\tilde{\Delta}^1_{\tau_z})^T \tilde{\Delta}^1_{\tau_z}.  $$
If $ \lim_{s \to \infty} E[(\tilde{\mathbf{x}}^k_s)^T \tilde{\mathbf{x}}^k_s ] < \infty $ then the above holds for $i = k-2$.
\end{thm}

\textbf{Proof of Theorem \ref{MainResult}:} \begin{enumerate}[(a)]
\item We know from Lemmas \ref{Markov} and \ref{Irreducible} that the process $\{(\mathbf{x}_s,\Delta_s)\}$ is an irreducible Markov chain. The set $S$ is small  (see Section \ref{sec:MCSC} and \cite{YukMeynTAC2010}). Using Lemma \ref{Drift} we can apply Theorem \ref{thm5} with $a=1$, the irreducible Markov chain $\{(\mathbf{x}_s, \Delta_s)\}$ and the functions $V(\mathbf{x}_s,\Delta_s)= (\Delta^1_s)^2,$ $\beta(\mathbf{x}_s,\Delta_s) =1$ and $b$ as given in Lemma \ref{Drift} to get that $\{(\mathbf{x}_s,\Delta_s)\}$ is positive Harris recurrent and has a unique invariant distribution.

\item Suppose that $\lambda \in \mathbb{R}$. We will apply Theorem \ref{thm5} with $a=0$, the irreducible Markov chain $\{(\mathbf{x}_s, \Delta_s)\}$ and the functions $V(\mathbf{x}_s,\Delta_s)= (\Delta^1_s)^2$, $ \beta(\mathbf{x}_s,\Delta_s) = \gamma (\Delta^1_s)^2$, $f(\mathbf{x}_s,\Delta_s)= \frac{\gamma}{\kappa} (x^n_s)^2$. From Lemma \ref{Drift}, we get
\begin{align*} &E[V(\mathbf{x}_{\tau_{z+1}},\Delta_{\tau_{z+1}}) \mid \mathcal{F}_{\tau_z}] = E[(\Delta^1_{\tau_{z+1}})^2 \mid \mathbf{x}_{\tau_z}, \Delta_{\tau_z}]\\
&\leq  (\Delta^1_{\tau_z})^2 -  \gamma E \left[ \left. \sum^{\tau_{z+1} -1}_{s=\tau_z} (\Delta^1_s)^2 \right \vert \mathbf{x}_{\tau_z}, \Delta_{\tau_z} \right] + b1_{ \{ (\mathbf{x}_{\tau_z}, \Delta_{\tau_z}) \in S \}} \nonumber \\
&\leq (\Delta^1_{\tau_z})^2 -  \gamma (\Delta^1_{\tau_z})^2 + b1_{ \{ (\mathbf{x}_{\tau_z}, \Delta_{\tau_z}) \in S \}} \\
&= V(\mathbf{x}_{\tau_z},\Delta_{\tau_z}) - \beta(\mathbf{x}_{\tau_z},\Delta_{\tau_z}) + b1_{ \{ (\mathbf{x}_{\tau_z}, \Delta_{\tau_z}) \in S \}}. \end{align*}
We know that Theorem \ref{MomentBound} holds immediately for $\{x^n_s\}$ and thus
\begin{align*} &E \left[ \left. \sum^{\tau_{z+1}-1}_{s=\tau_z} f(\mathbf{x}_s,\Delta_s) \right \vert \mathcal{F}_{\tau_z} \right] = \frac{\gamma}{\kappa}E \left[ \left. \sum^{\tau_{z+1}-1}_{s=\tau_z} (x^n_s)^2 \right \vert \mathbf{x}_{\tau_z}, \Delta_{\tau_z} \right] \leq \gamma (\Delta^1_{\tau_z})^2 = \beta(\mathbf{x}_{\tau_z},\Delta_{\tau_z}), \end{align*}
where we have used the ordering of bin sizes as described in Section \ref{sec:CCP}.

Thus, $ \lim_{s \to \infty}\frac{\gamma}{\kappa} E[(x^n_s)^2] < \infty $ by Theorem \ref{thm5} and so $ \lim_{s \to \infty} E[(x^n_s)^2] < \infty. $ This implies that Theorem \ref{MomentBound} holds for $\{x^{n-1}_s\}$ as mentioned in the proof and theorem statement. The finite second moment of all components then follows by induction.

In the complex case, we have that the drift condition \eqref{eq:drift} in Lemma \ref{Drift} also holds with $\Delta^2_s$ in place of $\Delta^1_s$ since they are equal. Choosing the functions $V(\mathbf{x}_s,\Delta_s)= (\Delta^1_s)^2 + (\Delta^2_s)^2, $ $ \beta(\mathbf{x}_s,\Delta_s)= \gamma ((\Delta^1_s)^2 + (\Delta^2_s)^2)$, $ f(\mathbf{x}_s,\Delta_s) = \frac{\gamma}{\kappa} (\tilde{\mathbf{x}}^n_s)^T\tilde{\mathbf{x}}^n_s$, we obtain the result. \end{enumerate}  \qed

\section{Multi-Sensor Systems} \label{sec:MSS}

\subsection{Problem Statement}

This is the main problem of the paper and is stated in Section \ref{sec:IntroPS}.

\subsection{Main Result} \label{sec:MultiMainResult}

To state the main result of this section, we first present a known result and an assumption.

The following theorem extends the classical observability canonical decomposition to the decentralized case. For a detailed proof in the centralized case, see \cite{Chen}. The more general multi-agent setup, where each agent makes observations and applies a control action, can be found in \cite{Morse}. We are not aware of an explicit proof and give a proof of Theorem \ref{BlockTheorem} in Section \ref{sec:PMSR} for the convenience of the reader.

\begin{thm} \label{BlockTheorem} Under Assumption \ref{ObservableAsn}, there exists a matrix $\mathbf{Q}$ such that if we define $\bar{\mathbf{A}} = \mathbf{Q} \mathbf{A} \mathbf{Q}^{-1}$ and $\bar{\mathbf{C}}^j = \mathbf{C}^j \mathbf{Q}^{-1}$ then
\begin{subequations} \label{eq:ABlock}
\begin{align} \bar{\mathbf{A}} =  \begin{bmatrix} \bar{\mathbf{A}}_M & * & \cdots & * \\
 & \bar{\mathbf{A}}_{M-1} & \cdots & * \\
 & & \ddots \\
 & 0 & & \bar{\mathbf{A}}_1  \end{bmatrix},  \label{eq:ABlock1} \\
\begin{bmatrix} \bar{\mathbf{C}}^M \\ \bar{\mathbf{C}}^{M-1} \\ \vdots \\ \bar{\mathbf{C}}^1 \end{bmatrix} =  \begin{bmatrix} \bar{\mathbf{C}}^M_O & * & \cdots & * \\
 & \bar{\mathbf{C}}^{M-1}_O & \cdots & * \\
 & & \ddots \\
 & 0 & & \bar{\mathbf{C}}^1_O  \end{bmatrix}, \label{eq:ABlock2} \end{align}
 \end{subequations}
 where the $*$'s denote irrelevant submatrices, each $\bar{\mathbf{A}}_j \in \mathbb{R}^{n_j \times n_j}$ and each $\bar{\mathbf{C}}^j_O \in \mathbb{R}^{p_j \times n_j}$.
\end{thm}

\begin{rem} \label{notuniqueorder} In the proof of Theorem \ref{BlockTheorem}, we give one construction for the triangular decomposition in \eqref{eq:ABlock}. This transformation is not unique. There may be many ways to achieve a block upper triangular form and it is not necessary to place the sensors in order $M,\dots,1$. \end{rem}


 Let us label the Jordan blocks of $\mathbf{A}$ as  $\mathbf{J}_1, \dots, \mathbf{J}_{\ell}$. Let $V_i$ be the (possibly generalized) eigenspace corresponding to $\mathbf{J}_i$. That is, if $\mathbf{v}_{i,1}, \dots, \mathbf{v}_{i,d_i}$ are the (possibly generalized) eigenvectors associated with $\mathbf{J}_i$ then $ V_i = \text{span} \{ \mathbf{v}_{i,1}, \dots, \mathbf{v}_{i,d_i} \}$ and has dimension $d_i$.

\begin{asn} \label{eigenspaces} Each eigenspace is observed by some sensor. That is, for each $1 \leq i \leq \ell$ there exists a $1 \leq j \leq M$ such that $V_i \subseteq O^j$. \end{asn}

The following is the main result of this section:

\begin{thm} \label{MultiSensorResult}
Under Assumption \ref{eigenspaces}, there exists a coding and control policy with average rate $R_{\text{avg}} \leq 1/(T2n) ( M + \sum^{n}_{i=1} \log_2( \lceil |\lambda_i|^{T2n} + \epsilon \rceil + 1) )$ for some $\epsilon > 0$ which gives:
\begin{enumerate}[(a)]
\item the existence of a unique invariant distribution for $\{\mathbf{x}_{2nt}\}$;
\item $ \lim_{t \to \infty} E[\|\mathbf{x}_{2nt}\|_2] < \infty. $
\end{enumerate}
\end{thm}

\begin{thm} \label{MultiMinRate} The average rate in Theorem \ref{MultiSensorResult} achieves the minimum rate \eqref{eq:rmin} asymptotically for large sampling periods. That is, $\lim_{T \to \infty} R_{\text{avg}} = R_{\min}$. \end{thm}

\textbf{Proof of Theorem \ref{MultiMinRate}: } Follows from the proof of Theorem \ref{MinimumRate}. \qed

\textbf{Proof of Theorem \ref{MultiSensorResult}:} Under Assumption \ref{eigenspaces}, we can assign each eigenspace $V_i \subseteq O^j$ to some sensor $j$. Let $V_{j,1}, \dots, V_{j,m_j}$ denote the eigenspaces assigned to sensor $j$ and let us write $V_{j,i} = \text{span}\{ \mathbf{v}_{j,i,1}, \dots, \mathbf{v}_{j,i,d_{j,i}} \}$ where each $\mathbf{v}_{j,i,h} \in \mathbb{R}^{n \times 1}$. We put
$ \mathbf{Q}_{j,i} = \begin{bmatrix} \mathbf{v}_{j,i,1} & \cdots & \mathbf{v}_{j,i,d_{j,i}} \end{bmatrix} ^T $, $\mathbf{Q}_j = \begin{bmatrix} \mathbf{Q}^T_{j,1} & \cdots & \mathbf{Q}^T_{j,m_j} \end{bmatrix}^T$ and $\mathbf{Q} = \begin{bmatrix} (\mathbf{Q}^M)^T & \cdots & (\mathbf{Q}^1)^T \end{bmatrix}^T$.

Each $\mathbf{v}_{j,i,h}$ belongs to the generalized eigenspace $V_{j,i}$, which is invariant under multiplication by $\mathbf{A}$. That is
\begin{align} \mathbf{v}_{j,i,h} \mathbf{A} \in V_{j,i}. \label{eq:invarianteigenspace} \end{align}
We apply the similarity transform $\bar{\mathbf{x}}_t = \mathbf{Q} \mathbf{x}_t$ to \eqref{eq:multisystem} and define $\bar{\mathbf{A}} = \mathbf{Q} \bar{\mathbf{A}} \mathbf{Q}^{-1}$, $\bar{\mathbf{B}} = \mathbf{Q} \mathbf{B}$ and $\bar{\mathbf{w}}_t = \mathbf{Q} \mathbf{w}_t$ to get the system
\begin{align} \bar{\mathbf{x}}_{t+1} = \bar{\mathbf{A}} \bar{\mathbf{x}}_t + \bar{\mathbf{B}} \bar{\mathbf{u}}_t + \bar{\mathbf{w}}_t. \end{align}
Furthermore, we can write $\bar{\mathbf{A}} = \text{diag} (\bar{\mathbf{A}}_M, \dots, \bar{\mathbf{A}}_1)$ where $\bar{\mathbf{A}}_j \in \mathbb{R}^{n_j \times n_j}$ and $n_j = \sum^{m_j}_{i=1} d_{j,i}$ is the sum of the dimensions of $V_{j,1}, \dots, V_{j,m_j}$. Equation \eqref{eq:invarianteigenspace} is analogous to \eqref{eq:linalg1} in the proof of Theorem \ref{BlockTheorem} and from this proof, we obtain the desired diagonal form.

We now look at the estimation of the state by the sensors. For convenience, let us write $\bar{\mathbf{x}}_t = \begin{bmatrix} (\bar{\mathbf{x}}^M_t)^T & \cdots & (\bar{\mathbf{x}}^1_t)^T\end{bmatrix}^T$ where $\bar{\mathbf{x}}^j_t = \begin{bmatrix} (\bar{\mathbf{x}}^{j,1}_t)^T & \cdots & (\bar{\mathbf{x}}^{j,m_j}_t)^T \end{bmatrix}^T $ and $ \bar{\mathbf{x}}^{j,i}_t = \begin{bmatrix} \bar{x}^{j,i,1}_t & \cdots & \bar{x}^{j,i,d_{j,i}}_t \end{bmatrix}^T $ with $ \bar{x}^{j,i,h}_t \in \mathbb{R}$. Let us write $ \mathcal{O}_{(C^j,A)} = \begin{bmatrix} (\mathbf{o}_{j,1})^T & \cdots & (\mathbf{o}_{j,np_j})^T \end{bmatrix}^T $ where each $\mathbf{o}_{j,i} \in \mathbb{R}^{1 \times n}$.

With our construction above, under Assumption \ref{eigenspaces}, we have for each $j,i,h$ that $ \mathbf{v}_{j,i,h} = \sum^{np_j}_{\ell=1} k^{j,i,h}_{\ell} \mathbf{o}_{j,\ell}$ for some real coefficients $\{k^{j,i,h}_{\ell}\}$. Consider the first $n$ time stages. By putting $\mathbf{k}^{j,i,h} = \begin{bmatrix} k^{j,i,h}_1 & \cdots & k^{j,i,h}_{np_j} \end{bmatrix}$, it follows that
\begin{align} &\mathbf{k}^{j,i,h} \begin{bmatrix} (\mathbf{y}^j_0)^T & \cdots & (\mathbf{y}^j_{n-1})^T \end{bmatrix}^T = \mathbf{k}^{j,i,h} \mathcal{O}_{(\mathbf{C}^j, \mathbf{A})} \mathbf{x}_0 + \bar{v}^{j,i,h}_0 \nonumber \\
&= \sum^{np_j}_{\ell=1} k^{j,i,h}_{\ell} \mathbf{o}_{j,\ell} \mathbf{x}_0  + \bar{v}^{j,i,h}_0 = \mathbf{v}_{j,i,h} \mathbf{x}_0 + \bar{v}^{j,i,h}_0 = \bar{x}^{j,i,h}_0 + \bar{v}^{j,i,h}_0 \nonumber \end{align}
where $\bar{v}^{j,i,h}_0$ is some zero mean Gaussian noise. We will use the same notation for $\bar{\mathbf{v}}^j_0$ that we use for $\bar{\mathbf{x}}^j_t$.

As in Lemma \ref{BeliefLemma} for the single-sensor case, we can use the next $n$ times stages to apply a control action. We then apply the above scheme repeatedly and sample every $2n$ time stages. By a slight abuse of notation, we define $\mathbf{x}_s = \bar{\mathbf{x}}_{2ns}$, $\mathbf{A} = \text{diag}(\mathbf{A}_M, \dots, \mathbf{A}_1) = \bar{\mathbf{A}}^{2n}$, $\mathbf{u}_s = \bar{\mathbf{u}}_{2ns}$, $\mathbf{w} = \bar{\mathbf{w}}_{2ns}$ and $\mathbf{v}^j_s = \bar{\mathbf{v}}^j_{2ns}$ to get the system
\begin{align} \mathbf{x}_{s+1} = \mathbf{A} \mathbf{x}_s + \mathbf{u}_s + \mathbf{w}_s, \quad \mathbf{y}^j_s = \mathbf{x}^j_s + \mathbf{v}^j_s, \quad 1 \leq j \leq M, \label{eq:diagmultisystem} \end{align}
where $\mathbf{x}_s = \begin{bmatrix} (\mathbf{x}^M_s)^T & \cdots & (\mathbf{x}^1_s)^T \end{bmatrix}^T$, $\mathbf{u}_s$ is chosen arbitrarily by the contoller and $\mathbf{y}^j_s$ is known by sensor $j$ at time $s$. The noise processes $\{\mathbf{w}_s\}$, $\{\mathbf{v}^j_s\}$ are each i.i.d. sequences of zero mean Gaussian random vectors. At time $s$, $\mathbf{w}_s$ and each $\mathbf{v}^j_s$ are independent of the state $\mathbf{x}_s$ but may be correlated with eachother. For $s_1 \neq s_2$ we have that $\mathbf{w}_{s_1}$ and $\mathbf{v}_{s_2}$ are independent.

Finally, we can assume that each $\mathbf{A}_j$ is in real Jordan form. Using the same notation for $\mathbf{x}_s, \Delta_s$ as for $\bar{\mathbf{x}}_t$, we associate with each $x^{j,i,h}_s$ the bin size $\Delta^{j,i,h}_s$. We define the sequence of stopping times
$$ \tau_0 = 0, \tau_{z+1} = \min\{s > \tau_z : |\mathbf{y}_s| = |\mathbf{x}_s + \mathbf{v}_s| \leq \Delta_s \}. $$
The feedback value $b_{2ns}$ is chosen as
$$ b_{2ns} = \begin{cases} 1, &\text{if $s=\tau_z$ for some $z \in \mathbb{N}$}, \\ 0, & \text{otherwise}, \end{cases} $$
so that we can then apply the same coding and control policy as in Section \ref{sec:CCP}. This reduces the problem to the single-sensor case and we obtain the result. \qed

\subsection{Sufficient Conditions for the General Multi-Sensor Case}

In Section \ref{sec:MultiMainResult}, Assumption \ref{eigenspaces} allowed us to diagonalize $\bar{\mathbf{A}}$ in \eqref{eq:ABlock} in Theorem \ref{BlockTheorem}. Without this assumption, the lower components of the state act as noise for the upper components. In particular, we need to bound these lower modes when all quantizers are perfectly-zoomed to achieve (b) of Theorem \ref{GeometricDecay}. To do this, we must have that the bin sizes of the lower modes are small compared with the upper ones. With many different eigenvalues, we cannot guarantee this in the general case. Below, we give a sufficient rate and an alternative assumption for stability.

For Theorem \ref{sufficientrate} below, let us write $\Lambda(\bar{\mathbf{A}}_j) = \{\lambda_{j,1}, \dots, \lambda_{j,n_j}\}$ where $\bar{\mathbf{A}}_j$ is given in \eqref{eq:ABlock}.

\begin{thm} \label{sufficientrate} There exists a coding and control policy which gives:
\begin{enumerate}[(a)]
\item the existence of a unique invariant distribution for $\{\mathbf{x}_{2nt}\}$;
\item $ \lim_{t \to \infty} E[\|\mathbf{x}_{2nt}\|_2] < \infty$,
\end{enumerate}
and with average rate in the limit of large sampling periods
$$\lim_{T \to \infty} R_{\text{avg}} = \sum^{M}_{j=1} \sum^{n_j}_{i=1} \log_2(\max \{|\lambda_{j,i}|, |\lambda_{h, \ell}| : h < j, 1 \leq \ell \leq n_h \}).$$
 \end{thm}
\textbf{Proof of Theorem \ref{sufficientrate}: } The proof follows that of Theorem \ref{MainResult}. The main difference is that we define $\lambda'_{j,i} = \max \{|\lambda_{j,i}|, |\lambda_{h, \ell}| : h < j, 1 \leq \ell \leq n_h \}$ and the bin numbers $K_{j,i} = \lceil (\lambda'_{j,i})^{2n} + \epsilon \rceil$ for some $\epsilon > 0$ and treat the lower components of the state as noise. \qed

Cleary, we could also achieve (a) and (b) in Theorem \ref{sufficientrate} with $\lim_{T \to \infty} R_{\text{avg}} = n \log_2(\lambda_{\text{absmax}})$ where $\lambda_{\text{absmax}} = \max_{j,i} \{|\lambda_{j,i}|\}$.

For Theorem \ref{orderedeigenvalues} below, recall that we have some flexibility in the decomposition given by Theorem \ref{BlockTheorem}. See the proof of Theorem \ref{BlockTheorem} and Remark \ref{notuniqueorder}.

\begin{thm} \label{orderedeigenvalues} If the eigenvalues of $\bar{\mathbf{A}}_M, \dots, \bar{\mathbf{A}}_1$ in \eqref{eq:ABlock} are ordered in decreasing magnitude then Theorem \ref{MultiSensorResult} holds without Assumption \ref{eigenspaces}. That is, the theorem holds if for $\lambda_i \in \Lambda(\bar{\mathbf{A}}_i)$ and $\lambda_j \in \Lambda(\bar{\mathbf{A}}_j)$ we have that $|\lambda_i| \leq |\lambda_j|$ when $i < j$.   \end{thm}

\textbf{Proof of Theorem \ref{orderedeigenvalues}: } The proof follows that of Theorem \ref{MainResult}. Since the eigenvalues are ordered in decreasing magnitude, we can maintain the ordering of the bin sizes and treat the lower components as noise. \qed

Finally, a remark on the vector scheme we have employed in this paper is in order.

\begin{rem} In this paper we present a vector stabilization scheme. From the problem statement, it would be natural to adopt a sequential stabilization scheme. That is, each of the components of the state is viewed as a separate system. In this case, we lose the Markov property and the number of time stages we must wait (denoted by $H$ in Theorem \ref{GeometricDecay}) to establish geometric decay is not uniform across the set of valid conditions $(\mathbf{x}_{\tau_z}, \Delta_{\tau_z})$. Such a scheme is left for future work.   \end{rem}

\section{Conclusion}

In this paper, we have presented a coding and control policy which achieves the minimum rate asymptotically in the limit of large sampling periods. We extend this result to the multi-sensor case under the assumption that each eigenspace is observed by some sensor. In the absence of this assumption, we give sufficient conditions for achieving stability. In all cases, we establish the existence of a unique invariant distribution for the sampled state and a finite second moment of the state. These strong forms of stability have not been considered in the literature for such systems to our knowledge. The proofs use random-time drift criteria for Markov chains. We wish to extend the results for control over general noisy channels along the lines of \cite{Yuksel2012}, \cite{Martins2}, \cite{Matveev}, \cite{SahaiParts} and \cite{Coviello2011}.

\section{Appendix}

\subsection{Supporting Results for Section \ref{sec:OPSSS}} \label{sec:SRSSS}

\textbf{Proof of Theorem \ref{MinimumRate}:} Let $\{\lambda_1, \dots, \lambda_n\}$ denote the set of eigenvalue of $\mathbf{A}$ with multiplicity. We give our control policy for period $T=2n$ with a fixed average rate of $R_{\text{avg}} = {\frac{1}{2n}}  \log_2 \left( \left\{ \prod^{n}_{i=1} \lceil |\lambda_i|^{2n} + \epsilon \rceil \right\} + 1 \right). $ Suppose that instead of sending an estimate every $2n$ time stages, we apply them periodically every $T2n$ time stages. Taking the limit as $T$ approaches infinity, our average rate satisfies
\begin{align*} &\lim_{T \to \infty} R_{avg} \leq \lim_{T \to \infty} \frac{1}{T2n} \left(  \sum^n_{i =1} \log_2(\lceil |\lambda_i|^{T2n} + \epsilon \rceil + 1) \right) \\
& = \lim_{T \to \infty}  \left(  \sum^n_{i =1} \log_2(\lceil |\lambda_i|^{T2n} + \epsilon \rceil + 1)^{\frac{1}{T2n}}  \right)  = \sum^{n}_{i=1} \log_2(|\lambda_i|).
 \end{align*}
In this sense, our policy achieves the minimum rate \eqref{eq:rmin} asymptotically. \qed

\textbf{Proof of Lemma \ref{BeliefLemma}:} Recall the basic recursion for LTI systems.
\begin{align*} \mathbf{x}_t &= \mathbf{A} \mathbf{x}_{t-1} + \mathbf{B} \mathbf{u}_{t-1} + \mathbf{w}_{t-1} = \mathbf{A}^2 \mathbf{x}_{t-2} + \mathbf{A}  \mathbf{B} \mathbf{u}_{t-2} + \mathbf{B} \mathbf{u}_{t-1} + \mathbf{A} \mathbf{w}_{t-2} + \mathbf{w}_{t-1} \\
& \dots = \mathbf{A}^t \mathbf{x}_0 + \sum^{t-1}_{i=0} \mathbf{A}^{t-1-i} \mathbf{B} \mathbf{u}_i + \sum^{t-1}_{i=0} \mathbf{A}^{t-1-i}\mathbf{w}_i.   \end{align*}
In the first $n$ time stages the sensor makes observations on the state and forms an estimate. In the second $n$ time stages we allow the controller to apply a control action.

We set $\mathbf{u}_i = 0$ for $0 \leq i \leq n-2$ so that the first $n$ observations of the sensor are
$$ \begin{bmatrix} \mathbf{y}_0 \\ \mathbf{y}_1 \\ \vdots \\ \mathbf{y}_{n-1} \end{bmatrix} = \mathcal{O}_{(\mathbf{C},\mathbf{A})}\mathbf{x}_0  + \begin{bmatrix} 0 \\ \mathbf{C} \mathbf{w}_0 \\ \vdots \\ \sum^{n-2}_{i=0} \mathbf{C} \mathbf{A}^{n-2-i}\mathbf{w}_i \end{bmatrix} +  \begin{bmatrix} \mathbf{v}_0 \\ \mathbf{v}_1 \\ \vdots \\ \mathbf{v}_{n-1} \end{bmatrix},$$
where $\mathcal{O}_{(\mathbf{C},\mathbf{A})}$ is the observability matrix of the pair $(\mathbf{C},\mathbf{A})$. We have assumed that $(\mathbf{C},\mathbf{A})$ is an observable pair. Equivalently, $\mathcal{O}_{(\mathbf{C},\mathbf{A})}$ has full column rank. By choosing a subset of $n$ equations from the matrix equation above, it is clear that we can apply the inverse to obtain the estimate $ \hat{\mathbf{y}}_0 = \mathbf{x}_0 + \sum^{n-2}_{i=0} \xi_i \mathbf{w}_i + \sum^{n-1}_{i=0} \zeta_i \mathbf{v}_i,  $ for some set $\{\xi_i,\zeta_i\}$ of matrices $\xi_i \in \mathbb{R}^{n \times n}$ and $\zeta_i \in \mathbb{R}^{n \times p}$.

Our estimate $\hat{\mathbf{y}}_0$ is generated at time $n-1$. At this time stage, the sensor sends the encoded value $q_{n-1}$ to the controller through the finite capacity channel. Based on this information, we allow the controller to apply control actions in time stages $n$ to $2n-1$. This is standard and we do not describe it in detail. We then have the system of equations
$$ \mathbf{x}_{2n} = \mathbf{A}^{2n} \mathbf{x}_{0} + \tilde{\mathbf{u}}_0 + \sum^{2n-1}_{i=0} \mathbf{A}^{2n-1-i}\mathbf{w}_i, \quad \hat{\mathbf{y}}_0 = \mathbf{x}_0 + \sum^{n-2}_{i=0} \xi_i \mathbf{w}_i + \sum^{n-1}_{i=0} \zeta_i \mathbf{v}_i, $$
where at time $n-1$, the estimate $\hat{\mathbf{y}}_0$ is known by the sensor and the action $\tilde{\mathbf{u}}_0$ is chosen arbitrarily by the controller.

Let us define the sampled variables $\tilde{\mathbf{x}}_s = \mathbf{x}_{2ns}$ and $\tilde{\mathbf{y}}_s = \hat{\mathbf{y}}_{2ns} $.
We define the noise processes
$$ \tilde{\mathbf{w}}_s =\sum^{2n-1}_{i=0} \mathbf{A}^{2n-1-i}\mathbf{w}_{2ns + i}, \quad \tilde{\mathbf{v}}_s = \sum^{n-2}_{i=0} \xi_i \mathbf{w}_{2ns + i} + \sum^{n-1}_{i=0} \zeta_i \mathbf{v}_{2ns + i}, $$
and note that they are both sequences of i.i.d. multivariate Gaussian random vectors. Then, by repeating our procedure every $2n$ time stages, we obtain the system
$$ \tilde{\mathbf{x}}_{s+1} = \mathbf{A}^{2n} \tilde{\mathbf{x}}_s + \tilde{\mathbf{u}}_s + \tilde{\mathbf{w}}_s,  \quad \tilde{\mathbf{y}}_{s} = \tilde{\mathbf{x}}_{s} + \tilde{\mathbf{v}}_s. $$
Finally, we apply a real Jordan transformation to the above system. We define $ \bar{\mathbf{x}}_s = \mathbf{P} \tilde{\mathbf{x}}_s$,   $\bar{\mathbf{A}} = \mathbf{P} \mathbf{A}^{2n} \mathbf{P}^{-1}$, $\bar{\mathbf{u}}_s = \mathbf{P} \tilde{\mathbf{u}}_s$, $\bar{\mathbf{w}}_s = \mathbf{P} \tilde{\mathbf{w}}_s$, $\bar{\mathbf{y}}_s = \mathbf{P}^{-1} \tilde{\mathbf{y}}_s$ and $\bar{\mathbf{v}}_s = \mathbf{P}^{-1} \tilde{\mathbf{v}}_s$ where $\mathbf{P}$ is the Jordan transform matrix. This gives the system
$$ \bar{\mathbf{x}}_{s+1} = \bar{\mathbf{A}} \bar{\mathbf{x}}_s + \bar{\mathbf{u}}_s + \bar{\mathbf{w}}_s, \quad \hat{\bar{\mathbf{x}}}_s = \bar{\mathbf{x}}_s + \bar{\mathbf{v}}_s. $$
Note that the matrix $\bar{\mathbf{A}}$ has eigenvalues $\lambda^{2n}_1,\dots,\lambda^{2n}_n$. \qed

\begin{rem} The estimate used in Lemma \ref{BeliefLemma} may appear naive. At first glance it would appear better to apply the Kalman filter. In this case, a new system is formed with the estimate as the state. The problem is that the noise for this system is not independent across time and we cannot extend our result to the multi-sensor case.
\end{rem}

\textbf{Proof of Lemma \ref{Markov}:}
Note that under our control policy we can write $ \mathbf{u}_s = g(\mathbf{x}_s, \mathbf{v}_s, \Delta_s)$ and $\Delta_{s+1} = f(\mathbf{x}_s, \mathbf{v}_s, \Delta_s) $ for some functions $g$ and $f$.

Let $\mathcal{B}(\mathbb{R}^n \times \mathbb{R}^n_+)$ be the Borel $\sigma$-field on $\mathbb{R}^n \times \mathbb{R}^n_+$.  It follows that
\begin{align*} & P \left(  (\mathbf{x}_{s+1},\Delta_{s+1}) \in (C \times D) \mid (\mathbf{x}_s,\Delta_s),\dots,(\mathbf{x}_0,\Delta_0) \right) \\
&= P \left(  \mathbf{x}_{s+1} \in C  \mid \Delta_{s+1} \in D, (\mathbf{x}_s,\Delta_s),\dots,(\mathbf{x}_0,\Delta_0) \right) \\
& \qquad P \left(  \Delta_{s+1} \in D \mid (\mathbf{x}_s,\Delta_s),\dots,(\mathbf{x}_0,\Delta_0) \right)  \end{align*}
\begin{align*}&= P \left(  \mathbf{A} \mathbf{x}_s + \mathbf{u}_s + \mathbf{w}_s \in C  \mid \Delta_{s+1} \in D, (\mathbf{x}_s,\Delta_s),\dots,(\mathbf{x}_0,\Delta_0) \right)  \\
&\qquad P \left(  f(\mathbf{x}_s, \mathbf{v}_s, \Delta_s) \in D \mid (\mathbf{x}_s,\Delta_s),\dots,(\mathbf{x}_0,\Delta_0) \right)  \end{align*}
\begin{align*}&= P \left(  \mathbf{A} \mathbf{x}_{s} + g(\mathbf{x}_s, \mathbf{v}_s, \Delta_s) + \mathbf{w}_s \in C  \mid \Delta_{s+1} \in D, (\mathbf{x}_s,\Delta_s),\dots,(\mathbf{x}_0,\Delta_0) \right)  \\
&\qquad P \left(  f(\mathbf{x}_s, \mathbf{v}_s, \Delta_s) \in D \mid (\mathbf{x}_s,\Delta_s),\dots,(\mathbf{x}_0,\Delta_0) \right)   \end{align*}
\begin{align*} &= P \left(   \mathbf{A}\mathbf{x}_s + g(\mathbf{x}_s, \mathbf{v}_s, \Delta_s) + \mathbf{w}_s \in C  \mid \Delta_{s+1} \in D, (\mathbf{x}_s,\Delta_s) \right) \\
& \qquad P \left(  f(\mathbf{x}_s, \mathbf{v}_s, \Delta_s) \in D \mid (\mathbf{x}_s,\Delta_s) \right)   \\
&= P \left(  (\mathbf{x}_{s+1},\Delta_{s+1}) \in (C \times D) \mid (\mathbf{x}_s,\Delta_s) \right), \end{align*}
for all $(C \times D ) \in \mathcal{B}(\mathbb{R}^n \times \mathbb{R}^n_+)$. \qed

\textbf{Proof of Lemma \ref{Irreducible}:} This follows immediately from the scalar case, as presented in the proof of Theorem 2.4 of \cite{Yuksel2010}. We can choose $\rho, \epsilon$ and $\eta$ such that $\log_2(\bar{Q}(q_s, \Delta_s))$ takes values in integer multiples of $\ell$ and the integers taken are relatively prime. By setting each $\Delta^i_0$ to be an integer multiple of $\ell$, it follows from the equation
$$ \log_2(\Delta^i_{s+1})/\ell = \log_2(\bar{Q}(q_s, \Delta_s))/\ell + \log_2(\Delta^i_s)/\ell $$
that $\log_2(\Delta^i_s)$ in an integer multiple of $\ell$ for all $s \in \mathbb{N}$. \qed

To prove Theorem \ref{GeometricDecay}, we need the following simple Gaussian bound. Recall that $\Lambda(\cdot)$ denotes the set of eigenvalues of its argument. Let us define $ \lambda_{\min} (\mathbf{A}) = \min \Lambda(\mathbf{A})$ and $\lambda_{\max} (\mathbf{A}) = \max \Lambda( \mathbf{A}).  $

\begin{lem} \label{SimpleGaussianBound}
Let $\mathbf{X} \sim \mathcal{N}(0,\mathbf{\Sigma})$ be a multivariate normal random variable with mean zero and covariance matrix $\mathbf{\Sigma} \in \mathbb{R}^{n \times n}$. For $\Delta \in \mathbb{R}^n_+$, the following bound holds.
$$ P(\mathbf{X} \nleq \Delta)
\leq 2  \sqrt{ \frac{ \lambda^{n+1}_{\max}(\mathbf{\Sigma})}{2\pi \det(\mathbf{\Sigma}) }} \sum^{n}_{i=1} \text{exp}\left\{-\frac{(\Delta^i)^2}{2\lambda_{\max}(\mathbf{\Sigma})}  \right\}. $$
\end{lem}

\textbf{Proof of Lemma \ref{SimpleGaussianBound}:}
Let $\mathbf{X} \sim \mathcal{N}(0,\mathbf{\Sigma})$ be a multivariate normal random vector with mean zero and covariance matrix $\mathbf{\Sigma} \in \mathbb{R}^{n \times n}$. We avoid the degenerate case and assume that $\mathbf{\Sigma}$ is positive-definite.
%
%
Let $\Delta \in \mathbb{R}^n_+$. Then
\begin{align} &P(\mathbf{X} \nleq \Delta) = P(\cup^{n}_{i=1}  \{ |x^i| > \Delta^i \}) \leq  \sum^{n}_{i=1} P(|x^i| > \Delta^i) \nonumber \end{align}
\begin{align}&= \sum^{n}_{i=1} \int_{|x^i| > \Delta^i} \frac{1}{\sqrt{(2\pi)^n \det(\mathbf{\Sigma})}} \text{exp}\left\{-\frac{1}{2} \mathbf{x}^T \mathbf{\Sigma}^{-1}\mathbf{x}\right\} \mathbf{dx} \nonumber \\
&\leq \sum^{n}_{i=1} \int_{|x^i| > \Delta^i} \frac{1}{\sqrt{(2\pi)^n \det(\mathbf{\Sigma})}} \text{exp}\left\{-\frac{1}{2} \lambda_{\min}(\mathbf{\Sigma}^{-1}) \mathbf{x}^T \mathbf{x}\right\} \mathbf{dx} \nonumber \end{align}
\begin{align} &=  \frac{1}{\sqrt{2\pi \det(\mathbf{\Sigma})\lambda^{n-1}_{\min}(\mathbf{\Sigma}^{-1})}} \sum^{n}_{i=1} 2 \int^{\infty}_{\Delta^i}  \text{exp}\left\{-\frac{1}{2} \lambda_{\min}(\mathbf{\Sigma}^{-1}) (x^i)^2\right\} dx^i \nonumber
\end{align}
\begin{align} &\leq \frac{2}{\sqrt{2\pi \det(\mathbf{\Sigma})\lambda^{n-1}_{\min}(\mathbf{\Sigma}^{-1})}} \sum^{n}_{i=1} \int^{\infty}_{\Delta^i} \frac{x^i}{\Delta_i} \text{exp}\left\{-\frac{1}{2} \lambda_{\min}(\mathbf{\Sigma}^{-1}) (x^i)^2\right\} dx^i \nonumber \\
&= \frac{2}{\sqrt{2\pi \det(\mathbf{\Sigma})\lambda^{n-1}_{\min}(\mathbf{\Sigma}^{-1})}} \sum^{n}_{i=1}  \frac{1}{\Delta^i} \left[- \frac{\text{exp}\left\{-\frac{1}{2} \lambda_{\min}(\mathbf{\Sigma}^{-1}) (x^i)^2\right\}}{\lambda_{\min}(\mathbf{\Sigma}^{-1})} \right]^{\infty}_{\Delta^i} \nonumber \end{align}
\begin{align} &= C \sum^{n}_{i=1}  \frac{1}{\Delta^i} \text{exp}\left\{-\frac{1}{2} \lambda_{\min}(\mathbf{\Sigma}^{-1}) (\Delta^i)^2 \right\} \leq C \sum^{n}_{i=1} \text{exp}\left\{-\frac{1}{2} \lambda_{\min}(\mathbf{\Sigma}^{-1}) (\Delta^i)^2 \right\}, \nonumber \end{align}
where the last line follows since we ensure $\Delta^i \geq 1$ for all $1 \leq i \leq n$ under our coding and control policy. We have also defined the constant $ C =   2/(\sqrt{2\pi \det(\mathbf{\Sigma}) \lambda^{n+1}_{\min}(\mathbf{\Sigma}^{-1})}). $  The eigenvalues of $\mathbf{\Sigma}^{-1}$ are the inverse eigenvalues of $\mathbf{\Sigma}$. This gives the desired bound.
\qed

\textbf{Proof of Theorem \ref{GeometricDecay}:} \textit{i)} \textbf{Exponential Bound.} Note that
\begin{align} &P(\tau_{z+1} - \tau_z > k \mid \mathbf{x}_{\tau_z}, \Delta_{\tau_z} ) = P\left( \left. \bigcap^{k}_{s=1} \left\{ |\mathbf{x}_{\tau_z + s} + \mathbf{v}_{\tau_z + s}| \nleq \frac{K}{2}\Delta_{\tau_z + s} \right \} \right \vert \mathbf{x}_{\tau_z}, \Delta_{\tau_z}   \right) \nonumber \\
&=  P\left( \left.  |\mathbf{x}_{\tau_z + k} + \mathbf{v}_{\tau_z + k}| \nleq \frac{K}{2}\Delta_{\tau_z + k}  \right \vert \bigcap^{k-1}_{s=1} \left \{ |\mathbf{x}_{\tau_z + s} + \mathbf{v}_{\tau_z + s}| \nleq \frac{K}{2}\Delta_{\tau_z + s} \right \}, \mathbf{x}_{\tau_z}, \Delta_{\tau_z} \right) \nonumber \end{align}
\begin{align} & \qquad \qquad P\left ( \left. \bigcap^{k-1}_{s=1} \left \{ |\mathbf{x}_{\tau_z + s} + \mathbf{v}_{\tau_z + s}| \nleq \frac{K}{2}\Delta_{\tau_z + s} \right \} \right \vert \mathbf{x}_{\tau_z}, \Delta_{\tau_z}  \right) \nonumber \\
&\leq P\left( \left.   |\mathbf{x}_{\tau_z + k} + \mathbf{v}_{\tau_z + k}| \nleq \frac{K}{2}\Delta_{\tau_z + k}   \right \vert \bigcap^{k-1}_{s=1} \left \{ |\mathbf{x}_{\tau_z + s} + \mathbf{v}_{\tau_z + s}| \nleq \frac{K}{2}\Delta_{\tau_z + s} \right \}, \mathbf{x}_{\tau_z}, \Delta_{\tau_z}\right) \nonumber \end{align}
\begin{align} &= P\left( \left. |\mathbf{x}_{\tau_z + k} + \mathbf{v}_{\tau_z + k}| \nleq \frac{K}{2}\Delta_{\tau_z + k} \right  \vert \tau_{z+1} - \tau_z > k-1, \mathbf{x}_{\tau_z}, \Delta_{\tau_z} \right). \label{eq:intbound1}  \end{align}
We first let $\lambda \in \mathbb{R}$. Let us define the noise vector $ \mathbf{w}_{\tau_z,k} = \frac{\mathbf{A}^k}{\lambda^k} ( - \mathbf{v}_{\tau_z} + \sum^{k-1}_{s=0} \mathbf{A}^{-1-s}\mathbf{w}_{\tau_z+s} ) + \frac{\mathbf{v}_{\tau_z + k}}{\lambda^k} $ and note that it is multivariate Gaussian. Before obtaining our bound, we define $ \xi =  \lceil |\lambda| + \epsilon \rceil / (|\lambda| + \epsilon - \eta) > 1. $ We let $N$ denote the nilpotent matrix (the matrix with ones on the upper diagonal) of appropriate size. Note that $N^{s} = 0$ for all $s \geq n$. Under our control policy, as described in Section \ref{sec:CCP}, we know that $ | (\mathbf{x}_{\tau_z} + \mathbf{v}_{\tau_z}) - \hat{\mathbf{x}}_{\tau_z}  | \leq \frac{1}{2} \Delta_{\tau_z}.$ It then follows that
\begin{align} &P\left( \left.  |\mathbf{x}_{\tau_z + k} + \mathbf{v}_{\tau_z + k}| \nleq \frac{K}{2}\Delta_{\tau_z + k}  \right \vert \tau_{z+1} - \tau_z > k-1,\mathbf{x}_{\tau_z}, \Delta_{\tau_z} \right) \nonumber \end{align}
\begin{align} &= P\Bigg(  \left \vert \mathbf{A}^k \mathbf{x}_{\tau_z} + \mathbf{A}^{k-1} \mathbf{u}_{\tau_z} + \sum^{k-1}_{s=0} \mathbf{A}^{k-1-s}\mathbf{w}_{\tau_z+s} + \mathbf{v}_{\tau_z + k} \right \vert \nonumber \\
& \qquad \qquad \nleq \frac{K}{2}\Delta_{\tau_z + k}  \Bigg \vert \tau_{z+1} - \tau_z > k-1, \mathbf{x}_{\tau_z}, \Delta_{\tau_z} \Bigg) \nonumber \\
&= P\Bigg(  \left \vert \mathbf{A}^k ( \mathbf{x}_{\tau_z} - \hat{\mathbf{x}}_{\tau_z} + \mathbf{v}_{\tau_z} - \mathbf{v}_{\tau_z} ) + \sum^{k-1}_{s=0} \mathbf{A}^{k-1-s}\mathbf{w}_{\tau_z+s} + \mathbf{v}_{\tau_z + k} \right \vert \nonumber \\
& \qquad \qquad  \nleq \frac{K}{2}\Delta_{\tau_z + k}  \Bigg \vert \tau_{z+1} - \tau_z > k-1, \mathbf{x}_{\tau_z}, \Delta_{\tau_z} \Bigg) \nonumber \\
&= P\Bigg( \left \vert (\lambda I + N)^k (\mathbf{x}_{\tau_z} + \mathbf{v}_{\tau_z} - \hat{\mathbf{x}}_{\tau_z} ) + \mathbf{A}^k \left( - \mathbf{v}_{\tau_z} + \sum^{k-1}_{s=0}  \mathbf{A}^{-1-s}\mathbf{w}_{\tau_z+s} \right ) + \mathbf{v}_{\tau_z + k} \right \vert \nonumber \\
& \qquad \qquad \nleq \frac{K}{2}\Delta_{\tau_z + k}  \Bigg \vert \tau_{z+1} - \tau_z > k-1, \mathbf{x}_{\tau_z}, \Delta_{\tau_z} \Bigg) \nonumber \end{align}
\begin{align} &\leq P\Bigg(  \left \vert \left\{  \lambda^k  +  \sum^{k}_{s=1} \binom{k}{s} \lambda^{k-s} N^{s}  \right \} (\mathbf{x}_{\tau_z} + \mathbf{v}_{\tau_z} - \hat{\mathbf{x}}_{\tau_z} ) +  \lambda^k\mathbf{w}_{\tau_z,k} \right \vert \nonumber \\
& \qquad \qquad \nleq \frac{K}{2} \rho^{k-1} |\lambda|^{k-1} \frac{|\lambda|}{|\lambda| + \epsilon - \eta} \Delta_{\tau_z} \Bigg \vert \mathbf{x}_{\tau_z}, \Delta_{\tau_z} \Bigg) \nonumber \end{align}
\begin{align} & \leq P\Bigg( \left\{ |\lambda|^k  +  \sum^{n}_{s=1} \binom{k}{s} |\lambda|^{k-s} N^{s} \right \}  |\mathbf{x}_{\tau_z} + \mathbf{v}_{\tau_z} - \hat{\mathbf{x}}_{\tau_z} | +  |\lambda|^k |\mathbf{w}_{\tau_z,k}| \nonumber \\
& \qquad \qquad \nleq \rho^{k-1} |\lambda|^{k} \xi \frac{1}{2} \Delta_{\tau_z}  \Bigg \vert \mathbf{x}_{\tau_z}, \Delta_{\tau_z} \Bigg) \nonumber \end{align}
\begin{align} & \leq P\Bigg( |\lambda|^k \left\{ \frac{1}{2} \Delta_{\tau_z}  + \delta \frac{1}{2} \Delta_{\tau_z} \sum^{n}_{s=1} \binom{k}{s} |\lambda|^{-s} +  |\mathbf{w}_{\tau_z,k}| \right \}  \nleq \rho^{k-1} |\lambda|^{k} \xi \frac{1}{2} \Delta_{\tau_z}  \Bigg \vert  \mathbf{x}_{\tau_z}, \Delta_{\tau_z} \Bigg) \label{eq:orderingforgeo} \\
& \leq P\Bigg(  \frac{1}{2} \Delta_{\tau_z}  + \delta \frac{1}{2} \Delta_{\tau_z} n k^n   +  |\mathbf{w}_{\tau_z,k}|  \nleq \rho^{k-1}  \xi \frac{1}{2} \Delta_{\tau_z}  \Bigg \vert  \mathbf{x}_{\tau_z}, \Delta_{\tau_z} \Bigg) \nonumber \end{align}
\begin{align} & \leq P\Bigg( |\mathbf{w}_{\tau_z,k}|  \nleq (\rho^{k-1}  \xi -1 - \delta n k^n) \frac{1}{2} \Delta_{\tau_z}  \Bigg \vert \mathbf{x}_{\tau_z}, \Delta_{\tau_z} \Bigg) \label{eq:geobound1} \end{align}
\begin{align} & \leq P\Bigg( |\mathbf{w}_{\tau_z,k}|  \nleq (\rho')^{k-1} \frac{1}{2} \Delta_{\tau_z}  \Bigg \vert \mathbf{x}_{\tau_z}, \Delta_{\tau_z} \Bigg)  \leq 2  \sqrt{ \frac{ \lambda^{n+1}_{\max}(\mathbf{\Sigma}_{\tau_z,k})}{2\pi \det(\mathbf{\Sigma}_{\tau_z,k}) }} \sum^{n}_{i=1} \text{exp}\left\{-\frac{(\rho')^{2(k-1)} (\Delta^i)^2}{8 \lambda_{\max}(\mathbf{\Sigma}_{\tau_z,k} )}   \right\}, \label{eq:geobound2}
\end{align}
where \eqref{eq:orderingforgeo} follows from our bin ordering. Equations \eqref{eq:geobound1} and \eqref{eq:geobound2} hold for all $k \geq H$ for some $H$ sufficiently large and in the special case of $k=1$. In the case $k=1$ we choose $\delta$ sufficiently small such that $ \xi -1 - \delta n > 0$. Equation \eqref{eq:geobound2} holds for some $1 < \rho' < \rho$ since we need only show that $ \rho^{k-1}  \xi -1 - \delta n k^n > (\rho')^{k-1} $
for sufficiently large $k$ and this follows since $\lim_{k \to \infty} \rho^{k-1}/(\rho')^{k-1} = \infty$ and $\lim_{k \to \infty} (-1- \delta n k^n)/(\rho')^{k-1} = 0$ by L'H$\hat{\text{o}}$pital's rule. In \eqref{eq:geobound2}, we have used Lemma \ref{SimpleGaussianBound} with the zero mean Gaussian vector $\mathbf{w}_{\tau_z,k}$ and denoted its covariance matrix by $\mathbf{\Sigma}_{\tau_z,k}$. From \eqref{eq:geobound2}, we can see that (b) of Theorem \ref{GeometricDecay} holds.

In order to bound \eqref{eq:geobound2} further, we define the covariance matrices $ \Sigma_{\mathbf{v}} = E[\mathbf{v}_{s} \mathbf{v}^T_{s}]$, $\Sigma_{\mathbf{v}, \mathbf{w}} = E[\mathbf{v}_{s} \mathbf{w}^T_{s}]$ and $\Sigma_{\mathbf{w}} = E[\mathbf{w}_{s} \mathbf{w}^T_{s}]$. Then
\begin{align*} &\Sigma_{\tau_z,k} = E[\mathbf{w}_{\tau_z,k}\mathbf{w}^T_{\tau_z,k}] = E\left\{ \frac{\mathbf{A}^k}{\lambda^k} \left( - \mathbf{v}_{\tau_z} + \sum^{k-1}_{s=0} \mathbf{A}^{-1-s}\mathbf{w}_{\tau_z+s} \right) + \frac{\mathbf{v}_{\tau_z + k}}{\lambda^k} \right\} \nonumber \\
& \qquad \left\{ \frac{\mathbf{A}^k}{\lambda^k} \left( - \mathbf{v}_{\tau_z} + \sum^{k-1}_{s=0} \mathbf{A}^{-1-s}\mathbf{w}_{\tau_z+s} \right) + \frac{\mathbf{v}_{\tau_z + k}}{\lambda^k} \right\}^T \\
&= \frac{\mathbf{A}^k}{\lambda^k} \left \{ \Sigma_{\mathbf{v}} - \Sigma_{\mathbf{v}, \mathbf{w}} (\mathbf{A}^{-1})^T -  \mathbf{A}^{-1}\Sigma^T_{\mathbf{v}, \mathbf{w}}+  \sum^{k-1}_{s=0} \mathbf{A}^{-1-s} \Sigma_{\mathbf{w}} (\mathbf{A}^{-1-s})^T \right \} \frac{(\mathbf{A}^k)^T}{\lambda^k} + \frac{\Sigma_{\mathbf{v}}}{\lambda^{2k}},  \end{align*}
where we have used the independence of $\{\mathbf{w}_s\}$, $\{\mathbf{v}_s\}$ across time and the independence of $\mathbf{v}_{s_1}$ and $\mathbf{w}_{s_2}$ for $s_1 \neq s_2$. Since both processes are zero mean, the cross terms are zero.

Recall that
$$ \mathbf{A}^k = \begin{bmatrix} \lambda^k & \binom{k}{1} \lambda^{k-1}  & \cdots & \binom{k}{n-1} \lambda^{k-n+1} \\
 & \lambda^k & \ddots &  \vdots \\
 & 0 & \ddots & \binom{k}{1} \lambda^{k-1} \\
 &  & & \lambda^k \end{bmatrix} $$
and let $Tr(\cdot)$ denote the trace of its argument. We get that
$$ Tr(\mathbf{A}^k (\mathbf{A}^k)^T) = \sum^{n-1}_{\ell=0} \sum^{\ell}_{s=0} \binom{k}{s}^2 \lambda^{2(k-s)} \leq n  \sum^{n-1}_{s=0} \binom{k}{s}^2 \lambda^{2(k-s)} \leq n \lambda^{2k} \sum^{n-1}_{s=0} \binom{k}{s}^2 \leq \lambda^{2k} n^2 k^{2n}.$$
Similarly, we can see that $Tr(\mathbf{A}^{k-1-s} (\mathbf{A}^{k-1-s} )^T) \leq \lambda^{2k} n^2 k^{2n}$ for all $0 \leq s \leq k-1$.

Define $ \mathbf{\Sigma}_1 = E[ (- \mathbf{v}_{\tau_z} + \mathbf{A}^{-1}\mathbf{w}_{\tau_z} )(- \mathbf{v}_{\tau_z} + \mathbf{A}^{-1}\mathbf{w}_{\tau_z})^T ] = \mathbf{\Sigma}_{\mathbf{v}} - \mathbf{\Sigma}_{\mathbf{v}, \mathbf{w}} (\mathbf{A}^{-1})^T -  \mathbf{A}^{-1}\mathbf{\Sigma}^T_{\mathbf{v}, \mathbf{w}}+ \mathbf{A}^{-1} \mathbf{\Sigma}_{\mathbf{w}} (\mathbf{A}^{-1})^T$. For symmetric matrices, every eigenvalue has an eigenvector. Thus, for $\mathbf{\Sigma}_{\tau_z,k}$, there exists a vector of unit length $\mathbf{e} \in \mathbb{R}^n$ such that
\begin{align} &\lambda_{\max}(\mathbf{\Sigma}_{\tau_z,k}) = \mathbf{e}^T \mathbf{\Sigma}_{\tau_z,k} \mathbf{e} =  \frac{1}{\lambda^{2k}} (\mathbf{e}^T \mathbf{A}^k) \mathbf{\Sigma}_1 ((\mathbf{A}^k)^T \mathbf{e}) \nonumber \\
& \qquad \qquad + \frac{1}{\lambda^{2k}} \sum^{k-1}_{s=1} (\mathbf{e}^T  \mathbf{A}^{k-1-s}) \mathbf{\Sigma}_{\mathbf{w}} ( (\mathbf{A}^{k-1-s})^T \mathbf{e}) + \frac{1}{\lambda^{2k}} \mathbf{e}^T \mathbf{\Sigma}_{\mathbf{v}} \mathbf{e} \nonumber \\
&\leq \frac{1}{\lambda^{2k}} \lambda_{\max}(\mathbf{\Sigma}_1) \mathbf{e}^T \mathbf{A}^k (\mathbf{A}^k)^T \mathbf{e} \nonumber \\
& \qquad \qquad + \frac{1}{\lambda^{2k}} \lambda_{\max}(\mathbf{\Sigma}_{\mathbf{w}}) \sum^{k-1}_{s=1} \mathbf{e}^T  \mathbf{A}^{k-1-s} (\mathbf{A}^{k-1-s})^T \mathbf{e} + \frac{1}{\lambda^{2k}} \lambda_{\max}(\mathbf{\Sigma}_{\mathbf{v}}) \mathbf{e}^T  \mathbf{e} \nonumber \end{align}
\begin{align} & \leq \frac{1}{\lambda^{2k}} \lambda_{\max}(\mathbf{\Sigma}_1) \lambda_{\max}(\mathbf{A}^k (\mathbf{A}^k)^T) \mathbf{e}^T  \mathbf{e} \nonumber \\
& \qquad \qquad + \frac{1}{\lambda^{2k}} \lambda_{\max}(\mathbf{\Sigma}_{\mathbf{w}}) \lambda_{\max}(\mathbf{A}^{k-1-s} (\mathbf{A}^{k-1-s})^T) \sum^{k-1}_{s=1} \mathbf{e}^T \mathbf{e} + \lambda_{\max}(\mathbf{\Sigma}_{\mathbf{v}}) \nonumber \end{align}
\begin{align} & \leq \frac{1}{\lambda^{2k}} \lambda_{\max}(\mathbf{\Sigma}_1) Tr(\mathbf{A}^k (\mathbf{A}^k)^T)  + k \frac{1}{\lambda^{2k}} \lambda_{\max}(\mathbf{\Sigma}_{\mathbf{w}}) Tr(\mathbf{A}^{k-1-s} (\mathbf{A}^{k-1-s})^T) + \lambda_{\max}(\mathbf{\Sigma}_{\mathbf{v}}) \nonumber \\
& \leq n^2 k^{2n+1}( \lambda_{\max}(\mathbf{\Sigma}_1) + \lambda_{\max}(\mathbf{\Sigma}_{\mathbf{w}}) + \lambda_{\max}(\mathbf{\Sigma}_{\mathbf{v}})). \nonumber
\end{align}
Recall Minkowski's Determinant Theorem (see for example \cite{Roberts73}). For nonnegative definite $n \times n$ matrices $\mathbf{A}$ and $\mathbf{B}$, it follows that $ \det(\mathbf{A} + \mathbf{B}) \geq \det(\mathbf{A}) + \det(\mathbf{B}). $

Using the above bound and the identity $\det(\mathbf{A}) = \det{\mathbf{A}^T}$ we get
\begin{align} &\det(\mathbf{\Sigma}_{\tau_z,k}) \geq \frac{1}{\lambda^{2k}} \det(\mathbf{A}^k)^2 \det \left(\mathbf{\Sigma}_1 + \sum^{k-1}_{s=1} \mathbf{A}^{-1-s} \mathbf{\Sigma}_{\mathbf{w}} (\mathbf{A}^{-1-s})^T \right) + \frac{1}{\lambda^{2k}} \det(\mathbf{\Sigma}_{\mathbf{v}}) \nonumber \\
& \geq  \det (\mathbf{\Sigma}_1) + \det \left( \sum^{k-1}_{s=1} \mathbf{A}^{-1-s} \mathbf{\Sigma}_{\mathbf{w}} (\mathbf{A}^{-1-s})^T \right) \geq \det (\mathbf{\Sigma}_1) + \sum^{k-1}_{s=1} \det( \mathbf{A}^{-1-s} \mathbf{\Sigma}_{\mathbf{w}} (\mathbf{A}^{-1-s})^T )  \nonumber \end{align}
\begin{align} &= \det (\mathbf{\Sigma}_1) + \det(\mathbf{\Sigma}_{\mathbf{w}}) \sum^{k-1}_{s=1} (\lambda^{-1-s})^{2n} \geq \det (\mathbf{\Sigma}_1). \nonumber\end{align}
Defining the constants $c_1 = n^2 ( \lambda_{\max}(\mathbf{\Sigma}_1) + \lambda_{\max}(\mathbf{\Sigma}_{\mathbf{w}}) + \lambda_{\max}(\mathbf{\Sigma}_{\mathbf{v}}))$ and $c_2 = \det(\mathbf{\Sigma}_1) $ we have obtained the bounds
\begin{align} &\lambda_{\max}(\mathbf{\Sigma}_1) \leq  c_1 k^{2n+1}, \quad \det(\mathbf{\Sigma}_{\tau_z,k}) \geq c_2.  \label{eq:lowerbound} \end{align}
Combining \eqref{eq:lowerbound} with \eqref{eq:geobound2} we get that
\begin{align} &P(\tau_{z+1} - \tau_z > k \mid \mathbf{x}_{\tau_z}, \Delta_{\tau_z} ) \leq 2  \sqrt{ \frac{ (c_1 k^{2n+1})^{n+1} }{2\pi c_2 } } \sum^{n}_{i=1} \text{exp}\left\{-\frac{(\rho')^{2(k-1)} (\Delta^i)^2}{8 c_1 k^{2n+1} }   \right\} \nonumber \\
& \qquad \leq C k^{n^2 + \frac{3}{2}n + \frac{1}{2}}  \sum^{n}_{i=1} \text{exp}\left\{-\frac{(\rho')^{2(k-1)} (\Delta^i)^2}{8 c_1 k^{2n+1} }   \right\} \label{eq:exponentialbound} \end{align}
where $C$ is the appropriate constant.

\textit{ii)} \textbf{Geometric Bound.} Note that in \eqref{eq:exponentialbound} we have a double exponential in $k$ since $ (\rho')^{2(k-1)} = e^{(k-1)2 \log(\rho')}. $ Let $a,b, c >0$ and recall that $ \lim_{k \to \infty} e^k/((a+b)k^{c+1}) = \infty $ by L'H$\hat{\text{o}}$pital's rule. This means that for all $L > 0$ there exists an $N$ such that $e^k/((a+b)k^{c+1}) > L$, for all $k \geq N$. Thus, $ e^k/k^c > L(a + b)k$, for all $k \geq N$. Then choosing $N$ large enough so that $(a + b)N > 1$ and subtracting $(a+b)k$ we get
$ e^k/k^c - (a + b)k > (L-1)(a + b)k > L-1$, for all $k \geq N$. Therefore, we have that $ \lim_{k \to \infty} e^k/k^c - (a + b)k = \infty. $ Since $\log(k) \leq k$ for $k \geq 1$, comparing with the above we find that $ \lim_{k \to \infty} e^k/k^c - a\log(k) - bk = \infty. $

Let $Q(k)$ be a polynomial of finite degree $m$. We can write $Q(k) = a_0 + a_1 k + \cdots + a_m k^m,$
for some coefficients $a_0,\dots,a_m \in \mathbb{R}$. Let $r > 1$ and consider
\begin{align} &\lim_{k \to \infty} r^k Q(k) \text{exp} \left\{ -\frac{e^k}{ k^c} \right\} = \lim_{k \to \infty} \sum^{m}_{i=0} a_i k^i r^k \text{exp} \left\{ -\frac{e^k}{k^c} \right\} \nonumber \\
&= \lim_{k \to \infty} \sum^{m}_{i=0} a_i \text{exp} \left\{ -\frac{e^k}{k^c} + i\log(k)  + \log(r) k  \right\} = 0. \label{eq:polybound1} \end{align}
Combining \eqref{eq:exponentialbound} and \eqref{eq:polybound1} gives the result for $\lambda \in \mathbb{R}$. The case of $\lambda \in \mathbb{C}$ is similar and we omit it. \qed

\textbf{Proof of Lemma \ref{Drift}:} We take $\lambda \in \mathbb{R}$. The proof for $\lambda \in \mathbb{C}$ is identical. We use the abbreviations
$P_z(k) = P(\tau_{z+1} - \tau_z = k \mid \mathbf{x}_{\tau_z}, \Delta_{\tau_z})$,  $P_z(k \mid Y) = P(\tau_{z+1} - \tau_z = k \mid Y, \mathbf{x}_{\tau_z}, \Delta_{\tau_z})$,
$\bar{P}_z(k) = P(\tau_{z+1} - \tau_z > k \mid \mathbf{x}_{\tau_z}, \Delta_{\tau_z})$ and $\bar{P}_z(k \mid Y) = P(\tau_{z+1} - \tau_z > k \mid Y, \mathbf{x}_{\tau_z}, \Delta_{\tau_z})$.


Put $r > \rho^2|\lambda|^2$. Using Theorem \ref{GeometricDecay}, we can bound the first term in \eqref{eq:drift} using the law of iterated expectations as follows
\begin{align} & E \left[ \left. \sum^{\tau_{z+1} -1}_{s=\tau_z} (\Delta^1_s)^2 \right \vert \mathbf{x}_{\tau_z}, \Delta_{\tau_z} \right] 
= \sum^{\infty}_{k=1} P_z(k) \sum^{k -1}_{s=0} E[ (\Delta^1_{\tau_z+s})^2  \vert \tau_{z+1} - \tau_z = k,  \Delta^1_{\tau_z}] \nonumber \end{align}
\begin{align} &\leq \sum^{\infty}_{k=1} P_z(k) \sum^{k -1}_{s=0} \rho^{2s} |\lambda |^{2s} (\Delta^1_{\tau_z})^2 \leq  (\Delta^1_{\tau_z})^2 \sum^{\infty}_{k=1}  k P_z(k) \rho^{2k}|\lambda |^{2k}   \nonumber \\
&\leq  (\Delta^1_{\tau_z})^2  \left(\sum^{H}_{k=1}  k P_z(k) \rho^{2k} |\lambda |^{2k}   +  \sum^{\infty}_{k=H+1}  \left( \frac{\rho^2|\lambda |^{2}}{r} \right)^k \right) = (\Delta^1_{\tau_z})^2 G_1. \label{eq:driftboundfirst} \end{align}
We have defined $ G_1 = \sum^{H}_{k=1}  k P_z(k) |\lambda |^{2k}   +  \sum^{\infty}_{k=H+1} ( \rho^2|\lambda |^2/r )^k < \infty. $
The series on the right converges since it is geometric. Similarly, we can bound the term $E[(\Delta^1_{\tau_{z+1}})^2 \mid \mathbf{x}_{\tau_z}, \Delta_{\tau_z}]$. Using the law of total expectation, we get
\begin{align} &E[(\Delta^1_{\tau_{z+1}})^2 \mid \mathbf{x}_{\tau_z}, \Delta_{\tau_z}] = P_z(1)E[(\Delta^1_{\tau_{z+1}})^2 \mid \tau_{z+1} - \tau_z = 1, \Delta_{\tau_z}] \nonumber \\
&\quad + \bar{P}_z(1)E[(\Delta^1_{\tau_{z+1}})^2 \mid \tau_{z+1} - \tau_z>1,\mathbf{x}_{\tau_z}, \Delta_{\tau_z}] = P_z(1) E[(\Delta^1_{\tau_{z+1}})^2 \mid \tau_{z+1} - \tau_z = 1, \Delta_{\tau_z}] \nonumber \end{align}
\begin{align} &\quad + \bar{P}_z(1)E[E[(\Delta^1_{\tau_{z+1}})^2 \mid \tau_{z+1}- \tau_z>1, \tau_{z+1}-\tau_z, \mathbf{x}_{\tau_z}, \Delta_{\tau_z}] \mid \tau_{z+1} - \tau_z > 1, \mathbf{x}_{\tau_z}, \Delta_{\tau_z}] \nonumber \\
&\leq P_z(1) E[(\Delta^1_{\tau_{z+1}})^2 \mid \tau_{z+1} - \tau_z = 1, \Delta_{\tau_z}] + \bar{P}_z(1) \sum^{\infty}_{k=2} P_z(k \mid \tau_{z+1} - \tau_z > 1) \rho^2|\lambda|^2 (\Delta^1_{\tau_z})^2 \nonumber \\
&=  P_z(1) E[(\Delta^1_{\tau_{z+1}})^2 \mid \tau_{z+1} - \tau_z = 1, \Delta_{\tau_z}] + \bar{P}_z(1) G_2 (\Delta^1_{\tau_z})^2  \label{eq:driftbound}
\end{align}
where we have defined $G_2 = \sum^{\infty}_{k=2} P_z(k \mid \tau_{z+1} - \tau_z > 1) \rho^2 |\lambda|^2 < \infty.  $ Convergence comes from the geometric decay, as in the previous bound. Note that the geometric bound in Theorem \ref{GeometricDecay} still holds with $P_z(k \mid \tau_{z+1} - \tau_z > 1)$ in place of $P_z(k)$ since we obtain our bound by looking only at the $\tau_z + k$ term, as can be seen in \eqref{eq:intbound1}.

There exists a $\zeta$ such that $0 < \zeta < 1 - ( |\lambda|/(|\lambda| + \epsilon - \eta) )^2. $ We know from Theorem \ref{GeometricDecay} that $ \lim_{\Delta_{\tau_z} \to \infty} \bar{P}_z(1) = 0. $ Recall that $\Delta^i_s \geq \bar{L}^i$ for all $t \in \mathbb{N}$. Then, we choose $\mathbf{L}$ large enough to get an appropriate $\bar{\mathbf{L}}$ such that
$ \bar{P}_z(1) G_2 < \zeta. $ We put
$$ \gamma =  \frac{1 - \left(\frac{|\lambda|}{|\lambda| + \epsilon - \eta}\right)^2 - \zeta}{G_1} $$
so that $\gamma > 0$. Now, if $\Delta_{\tau_z} \notin S_{\Delta}$ then we have that $\Delta^1_{\tau_z} \geq F$ since $\Delta^1_s \geq \Delta^2_s \geq \cdots \geq \Delta^n_s$ for all $t \in \mathbb{N}$ by construction. Since $F > L^1$, the bin size shrinks and
$$ E[(\Delta^1_{\tau_{z+1}})^2 \mid \tau_{z+1} - \tau_z = 1, \Delta_{\tau_z}] = \left( \frac{|\lambda|}{|\lambda| + \epsilon - \eta}\right)^2 (\Delta^1_{\tau_z})^2. $$
If $\Delta_{\tau_z} \in S_{\Delta}$ then we use the simple bound $ E[(\Delta^1_{\tau_{z+1}})^2 \mid \tau_{z+1} - \tau_z = 1, \Delta_{\tau_z}] \leq  \rho^2|\lambda|^2 (\Delta^1_{\tau_z})^2. $

From the above, we have the following bounds. If $\Delta_{\tau_z} \notin S_{\Delta}$ then
\begin{align} E[(\Delta^1_{\tau_{z+1}})^2 \mid \mathbf{x}_{\tau_z}, \Delta_{\tau_z}] \leq (\Delta^1_{\tau_z})^2 \left\{ \left( \frac{|\lambda|}{|\lambda| + \epsilon - \eta} \right)^2 + \zeta \right\}. \label{eq:driftbound1}  \end{align}
If $\Delta_{\tau_z} \in S_{\Delta}$ then
\begin{align} E[(\Delta^1_{\tau_{z+1}})^2 \mid \mathbf{x}_{\tau_z}, \Delta_{\tau_z}] \leq (\Delta^1_{\tau_z})^2 \{ \rho^2 |\lambda|^2 + \zeta \}. \label{eq:driftbound2}  \end{align}
In the case $\Delta_{\tau_z} \notin S_{\Delta}$ we apply \eqref{eq:driftboundfirst}, \eqref{eq:driftbound} and \eqref{eq:driftbound1} to get
\begin{align} &\gamma E \left[ \left. \sum^{\tau_{z+1} -1}_{s=\tau_z} (\Delta^1_s)^2 \right \vert \mathbf{x}_{\tau_z}, \Delta_{\tau_z} \right] \leq (\Delta^1_{\tau_z})^2 \gamma G_1 \nonumber \\
&= (\Delta^1_{\tau_z})^2 \left\{ 1 - \left( \frac{|\lambda|}{|\lambda| + \epsilon - \eta} \right)^2 - \zeta \right\} \leq  (\Delta^1_{\tau_z})^2 - E[(\Delta^1_{\tau_{z+1}})^2|\mathbf{x}_{\tau_z}, \Delta_{\tau_z}]. \nonumber \end{align}
In the case $\Delta_{\tau_z} \in S_{\Delta}$ we apply \eqref{eq:driftboundfirst}, \eqref{eq:driftbound} and \eqref{eq:driftbound2} to get
\begin{align} &\gamma E \left[ \left. \sum^{\tau_{z+1} -1}_{s=\tau_z} (\Delta^1_s)^2 \right \vert \mathbf{x}_{\tau_z}, \Delta_{\tau_z} \right] \leq (\Delta^1_{\tau_z})^2 \gamma G_1 = (\Delta^1_{\tau_z})^2 \left\{ 1 - \left( \frac{|\lambda|}{|\lambda| + \epsilon - \eta}\right)^2 - \zeta \right\} \nonumber \\
&= (\Delta^1_{\tau_z})^2 - (\Delta^1_{\tau_z})^2\{\rho^2|\lambda|^2 + \zeta\}  + (\Delta^1_{\tau_z})^2 \left\{ \rho^2 |\lambda|^2 - \left(\frac{|\lambda|}{|\lambda| + \epsilon - \eta}\right)^2 \right\} \nonumber \\
&\leq (\Delta^1_{\tau_z})^2 - E[(\Delta^1_{\tau_{z+1}})^2 \mid \mathbf{x}_{\tau_z}, \Delta_{\tau_z}] + F^2 \left\{ \rho^2 |\lambda|^2 - \left(\frac{ |\lambda|}{|\lambda| + \epsilon - \eta}\right)^2 \right\}.  \nonumber \end{align}
We set $ b= F^2 \{ \rho^2|\lambda|^2 - ( |\lambda|/(|\lambda| + \epsilon - \eta) )^2 \}. $ Since $\Delta_{\tau_z} \in S_{\Delta}$ if and only if $(\mathbf{x}_{\tau_z}, \Delta_{\tau_z}) \in S$, we obtain Lemma \ref{Drift}. \qed


\textbf{Proof of Theorem \ref{MomentBound}:} Consider first the case $\lambda \in \mathbb{R}$ and let $x^{n+1}_s=0$. Using the law of total expectation we get
\begin{align} &E \left[ \left. \sum^{\tau_{z+1}-1}_{s= \tau_z} (x^i_s)^2 \right \vert \mathbf{x}_{\tau_z},\Delta_{\tau_z} \right] = E\Big[ E \Big[  (x^i_{\tau_z})^2 + \sum^{\tau_{z+1}-1}_{s= \tau_z + 1} (\lambda x^i_{s-1} + x^{i+1}_{s-1} + u^i_{s-1} + w^i_{s-1})^2 \nonumber \\
& \qquad \Big \vert \tau_{z+1} - \tau_z, \mathbf{x}_{\tau_z},\Delta_{\tau_z}  \Big] \Big \vert \mathbf{x}_{\tau_z},\Delta_{\tau_z} \Big] \nonumber \\
&= \sum^{\infty}_{k=1} P_z(k) E \Bigg[ (x^i_{\tau_z})^2 + \sum^{k-1}_{s=1} \Bigg(\lambda^s x^i_{\tau_z} + \lambda^{s-1} x^{i+1}_{\tau_z} - \lambda^{s-1} u^i_{\tau_z}  + \sum^{s-1}_{j=1} \lambda^{s-1-j} x^{i+1}_{\tau_z + j} \nonumber \\
& \qquad + \sum^{s-1}_{j=0} \lambda^{s-1-j} w^i_{\tau_z + j}\Bigg)^2 \Bigg \vert \mathbf{x}_{\tau_z},\Delta_{\tau_z} \Bigg] \nonumber \end{align}
\begin{align} &= \sum^{\infty}_{k=1} P_z(k) E \Bigg[ (x^i_{\tau_z})^2 + \sum^{k-1}_{s=1} \Bigg(\lambda^s (x^i_{\tau_z} + v^i_{\tau_z} -  \hat{x}^i_{\tau_z}) + \lambda^{s-1} (x^{i+1}_{\tau_{z}} + v^{i+1}_{\tau_z} - \hat{x}^{i+1}_{\tau_{z}})   \nonumber \\
& \qquad + \sum^{s-1}_{j=1} \lambda^{s-1-j} x^{i+1}_{\tau_z + j} - \lambda^s v^i_{\tau_z}- \lambda^{s-1} v^{i+1}_{\tau_z} + \sum^{s-1}_{j=0} \lambda^{s-1-j} w^i_{\tau_z + j}\Bigg)^2 \Bigg \vert \mathbf{x}_{\tau_z},\Delta_{\tau_z} \Bigg] \nonumber
\end{align}
\begin{align}
&\leq 6 \sum^{\infty}_{k=1} P_z(k) E \Bigg[ (x^i_{\tau_z})^2 + \sum^{k-1}_{s=1} \Big ( \lambda^{2s} (x^i_{\tau_z} + v^i_{\tau_z} -  \hat{x}^i_{\tau_z})^2 + \lambda^{2(s-1)} (x^{i+1}_{\tau_{z}} + v^{i+1}_{\tau_z} - \hat{x}^{i+1}_{\tau_{z}})^2   \nonumber \\
& \qquad + \left( \sum^{s-1}_{j=1} \lambda^{s-1-j} x^{i+1}_{\tau_z + j} \right)^{2} + \lambda^{2s} (v^i_{\tau_z})^2 + \lambda^{2(s-1)} (v^{i+1}_{\tau_z})^2 \nonumber \\
& \qquad + \left(\sum^{s-1}_{j=0} \lambda^{s-1-j} w^i_{\tau_z + j}\right)^2 \Big ) \Bigg \vert \mathbf{x}_{\tau_z},\Delta_{\tau_z} \Bigg] \label{eq:moment1} \end{align}
\begin{align} &\leq 6 \sum^{\infty}_{k=1} P_z(k) E \Bigg[ \frac{K^2}{4}(\Delta^i_{\tau_z})^2 + \sum^{k-1}_{s=1} \Big ( \lambda^{2s} \left(\frac{1}{2}\Delta^i_{\tau_z}\right)^2 + \lambda^{2(s-1)} \left(\frac{1}{2}\Delta^{i+1}_{\tau_z}\right)^2 \nonumber \\
& \qquad + s\sum^{s-1}_{j=1} \lambda^{2(s-1-j)} (x^{i+1}_{\tau_z + j})^2 + \lambda^{2s} (v^i_{\tau_z})^2 + \lambda^{2(s-1)} (v^{i+1}_{\tau_z})^2 \nonumber \\
& \qquad + s\sum^{s-1}_{j=0} \lambda^{2(s-1-j)} (w^i_{\tau_z + j})^2 \Big ) \Bigg \vert \mathbf{x}_{\tau_z},\Delta_{\tau_z} \Bigg] \label{eq:moment2} \end{align}
\begin{align} & \leq 6 \sum^{\infty}_{k=1} P_z(k) \Bigg ( \frac{K^2}{4} (\Delta^1_{\tau_z})^2 + \sum^{k-1}_{s=1} \Big ( \lambda^{2s} \frac{1}{2} \left(\Delta^1_{\tau_z}\right)^2 + s^2 \lambda^{2s} M^{i+1} \nonumber \\
& \qquad + \lambda^{2s} \sigma^2_{v,i} + \lambda^{2s} \sigma^2_{v,i+1} + s^2 \lambda^{2s} \sigma^2_{w,i} \Big ) \Bigg )   \nonumber \\
& \leq    (\Delta^1_{\tau_z})^2 6 \sum^{\infty}_{k=1} P_z(k) \left \{ \frac{K^2}{4} + \sum^{k-1}_{s=1}  \lambda^{2s} \left( \frac{1}{2} +  s^2  M^{i+1} + \sigma^2_{v,i} +  \sigma^2_{v,i+1} + s^2  \sigma^2_{w,i} \right) \right \}.  \label{eq:moment3}
\end{align}
In \eqref{eq:moment1} and \eqref{eq:moment2} we have used Jensen's inequality. Line \eqref{eq:moment3} follows since we can bound $\Delta^i_{s} > 1$ for all $s \in \mathbb{N}$. We have defined $ M^i = \sup_{s \in \mathbb{N}} E[(x^i_s)^2] < \infty$, $\sigma^2_{v,i} = E[(v^i_s)^2]$ and $\sigma^2_{w,i} = E[(w^i_s)^2]$. The fact that $M^{i}$ is finite for $2 \leq i \leq n-1$ follows from induction in the proof of Theorem \ref{MainResult}. By convention we put $M^{n+1}=0$. Now, we apply Theorem \ref{GeometricDecay} with $Q(k)=k^3$ and $r > \lambda_2$ to yield
\begin{align*} &\sum^{\infty}_{k=1} P_z(k) \sum^{k-1}_{s=1} \lambda^{2s} s^2 = \sum^{H}_{k=1}  \sum^{k-1}_{s=1} P_z(k) \lambda^{2s} s^2       + \sum^{\infty}_{k=H+1} \sum^{k-1}_{s=1} \lambda^{2s} P_z(k) s^2 \\
&\leq  G + \sum^{\infty}_{k=H+1}  \lambda^{2k} P_z(k) k^3 \leq  G + \sum^{\infty}_{k=H+1}\left( \frac{ \lambda^{2}}{r} \right)^k < \infty. \end{align*}
The last series converges since it is geometric. We have defined $G =\sum^{H}_{k=1}  \sum^{k-1}_{s=1} P_z(k) \lambda^{2s} s^2 < \infty.$ Therefore we can set
$$ \kappa = 6 \sum^{\infty}_{k=1} P_z(k) ( K^2/4 + \sum^{k-1}_{s=0}   \lambda^{2s} \left( 1/2 +  s^2  M^{i+1} + \sigma^2_{v,i} +  \sigma^2_{v,i+1} + s^2  \sigma^2_{w,i} \right) < \infty $$
to get the result.
For $\lambda \in \mathbb{C}$, the proof is similar and we omit it. \qed

\subsection{A Supporting Result for Section \ref{sec:MultiMainResult}} \label{sec:PMSR}

\textbf{Proof of Theorem \ref{BlockTheorem}:} We define $n_1 = \dim (O^1)$, and $n_j = \dim( O^j \backslash \{ \cup^{j-1}_{i=1} O^i \})$, for $2 \leq j \leq M$. We choose $n_1$ linearly independent row vectors from $\mathcal{O}_{(\mathbf{C}^1,\mathbf{A})}$ and label them $\mathbf{q}^1_1,\dots, \mathbf{q}^1_{n_1}$. 
Proceeding by induction, we choose $\{\mathbf{q}^j_1,\dots, \mathbf{q}^j_{n_j} \}$  from $\mathcal{O}_{(\mathbf{C}^j,\mathbf{A})}$ such that $$\{\mathbf{q}^1_1,\dots, \mathbf{q}^1_{n_1}, \mathbf{q}^2_1,\dots, \mathbf{q}^2_{n_2}, \dots, \mathbf{q}^j_1,\dots, \mathbf{q}^j_{n_j} \}$$
is a set of linearly independent vectors.

We define $ \mathbf{Q}^j = \begin{bmatrix} (\mathbf{q}^j_1)^T  & \cdots & (\mathbf{q}^j_{n_j})^T  \end{bmatrix}^T$  for all $1 \leq j \leq M$  and concatenate these matrices to choose our transformation matrix $ \mathbf{Q} = \begin{bmatrix} (\mathbf{Q}^M)^T & \cdots & (\mathbf{Q}^1)^T \end{bmatrix}^T. $ It will also be convenient to denote the rows of $\mathbf{Q}$ by $\mathbf{q}_1,\dots,\mathbf{q}_n$ so that $ \mathbf{Q} = \begin{bmatrix} (\mathbf{q}_1)^T & \cdots & (\mathbf{q}_n)^T \end{bmatrix}^T. $

From the Cayley-Hamilton Theorem, we know that for all $m \geq n$ there exist $\alpha_0,\dots,\alpha_{n-1}$ such that $ \mathbf{A}^m = \sum^{n-1}_{i=0} \alpha_i \mathbf{A}^i. $ Since $\{\mathbf{q}^j_1,\dots, \mathbf{q}^j_{n_j} \}$ are rows of $\mathcal{O}_{(\mathbf{C}^j,\mathbf{A})}$, this implies that $\mathbf{q}^j_i \mathbf{A}$ is in the row space of $\mathcal{O}_{(\mathbf{C}^j,\mathbf{A})}$ for all $1 \leq i \leq n_j$. Let us define the sets
$$ S_j := \{\mathbf{q}^1_1,\dots, \mathbf{q}^1_{n_1}, \mathbf{q}^2_1,\dots, \mathbf{q}^2_{n_2}, \dots, \mathbf{q}^j_1,\dots, \mathbf{q}^j_{n_j} \}, \qquad 1 \leq j \leq M$$
From our construction, it is then clear that
\begin{align} \mathbf{q}^j_1 \mathbf{A},\dots, \mathbf{q}^j_{n_j} \mathbf{A} \in \text{span} (S_j). \label{eq:linalg1} \end{align}
We write $\bar{\mathbf{A}}$ in terms of its column vectors as $ \bar{\mathbf{A}} = \begin{bmatrix} \bar{\mathbf{a}}_1 & \cdots & \bar{\mathbf{a}}_n \end{bmatrix}$ where $\bar{\mathbf{a}}_i \in \mathbb{R}^{n \times 1}$ for each $1 \leq i \leq n$. Our similarity transform gives
\begin{align} \bar{\mathbf{A}} \mathbf{Q} = \mathbf{Q} \mathbf{A}. \label{eq:linalg} \end{align}
Recall from linear algebra that we can write the left side of \eqref{eq:linalg} as $ \sum^{n}_{i=1} \bar{\mathbf{a}}_i \mathbf{q}_i$
where $\bar{\mathbf{a}}_i \mathbf{q}_i \in \mathbb{R}^{n \times n}$ for each $1 \leq i \leq n$. Now, to return to our earlier notation, each vector $\mathbf{q}^j_i \mathbf{A}$ is a linear combination of $\{\mathbf{q}^k_{\ell} : 1 \leq k \leq j, 1 \leq \ell \leq n_k  \} $ and is linearly independent of the remaining rows of $\mathbf{Q}$. Since $\mathbf{Q} \mathbf{A} = \begin{bmatrix} (\mathbf{q}_1 \mathbf{A})^T & \cdots & (\mathbf{q}_n \mathbf{A})^T \end{bmatrix}^T$, we see from \eqref{eq:linalg} that the $i^{th}$ row of $\bar{\mathbf{A}}$ is the representation of $\mathbf{q}_i \mathbf{A}$ with respect to $\mathbf{q}_1, \dots, \mathbf{q}_n$. More precisely, we write $\bar{\mathbf{a}}_i = \begin{bmatrix} \bar{a}_{i,1} & \cdots & \bar{a}_{i,n} \end{bmatrix}^T$ with each $\bar{a}_{i,h} \in \mathbb{R}$ so that \eqref{eq:linalg} gives the system of equations
\begin{align} \sum^{n}_{i=1} \bar{a}_{i,h} \mathbf{q}_i = \mathbf{q}_h \mathbf{A}, \quad 1 \leq h \leq n. \label{eq:linalg2} \end{align}
Combing \eqref{eq:linalg1} and \eqref{eq:linalg2} gives the desired form.

We next turn our attention to the form of $\bar{\mathbf{C}}^j$. Since each $\mathbf{C}^j$ is a submatrix of $\mathcal{O}_{(\mathbf{C}^j,\mathbf{A})}$, it is clear that the rows of $\mathbf{C}^j$ are in the span of $S_j$. Since $ \bar{\mathbf{C}}^j \mathbf{Q} = \mathbf{C}^j$, by writing $\bar{\mathbf{C}}^j$ in terms of its column vectors we obtain the desired form. \qed

\subsection{Review of Symmetric Matrices, Markov Chains and Stochastic Stability} \label{sec:MCSC}




Recall that a matrix $\mathbf{A} \in \mathbb{R}^{n \times n}$ is said to be \textit{symmetric} if $ \mathbf{A}^T = \mathbf{A}.$


\begin{lem} \label{MatrixLemma2} Let $\mathbf{A} \in \mathbb{R}^{n \times n}$ be a symmetric matrix with eigenvalues $\lambda_1,\dots,\lambda_n$. If we let $ \lambda_{\min} = \min \{ \lambda_1 ,\dots , \lambda_n\}$ and $ \lambda_{\max} = \max \{ \lambda_1 ,\dots , \lambda_n\} $ then $ \lambda_{\min} \mathbf{x}^T \mathbf{x} \leq \mathbf{x}^T  \mathbf{A} \mathbf{x} \leq \lambda_{\max} \mathbf{x}^T \mathbf{x} $ for all $\mathbf{x} \in \mathbb{R}^n$. \end{lem}

We present a brief discussion on stochastic stability of Markov chains. For a list of definitions on Markov chains, the reader is referred to \cite{MeynBook} and \cite{Meyn}. Let $\phi = \{ \phi_t, t\geq 0\}$ be a Markov chain with a complete separable metric state space $(\mathbb{X},{\cal B}(\mathbb{X}))$, and defined on a probability space $(\Omega, {\cal F}, {\cal P})$, where
 ${\cal B}(\mathbb{X})$ denotes the Borel $\sigma-$field on $\mathbb{X}$, $\Omega$ is the sample space, $ {\cal F}$ a sigma field of subsets of $\Omega$, and ${\cal P}$ a probability measure. Let $P(x,D) := P(\phi_{t+1}\in D | \phi_t=x)$ denote the transition probability from $x$ to $D$.

\begin{defn}
For a Markov chain, a probability measure $\pi$ is {\it invariant} on the Borel space $(\mathbb{X}, {\cal B}(\mathbb{X}))$ if $\pi(D)=\int_{\mathbb{X}} P(x,D)\pi(dx), \quad \forall D \in {\cal B}(\mathbb{X}).$
\end{defn}

\begin{defn}
A Markov chain is {\it $\mu$-irreducible} if for any set $B \in {\cal B}(\mathbb{X})$ with $\mu(B) > 0$ and $\forall x \in \mathbb{X}$,
there exists some integer $n>0$, possibly depending on $B$ and $x$, such that $P^n(x,B) > 0$, where $P^n(x,B)$ is the transition probability in $n$ stages. That is $P(\phi_{t+n} \in B | \phi_t=x)$.
\end{defn}

\begin{defn}
A set $A \subset \mathbb{X}$ is small if there is an integer $n \geq 1$ and a positive measure $\mu$ satisfying $\mu(\mathbb{X}) > 0$ and $P^n(x,B) \geq \mu(B), \quad \forall x \in A, B \in {\cal B}(\mathbb{X}).$
\end{defn}

\begin{defn}
A set $A \subset \mathbb{X}$ is $\zeta-$petite on $(\mathbb{X},{\cal B}(\mathbb{X}))$ if for
 some distribution ${\cal Z}$ on $\mathbb{N}$ (set of natural numbers), and some non-trivial measure $\zeta$,
$\sum^{\infty}_{n=0}P^n(x,B){\cal Z}(n) \geq \zeta(B)$, $\forall x \in A, B \in {\cal B}(\mathbb{X}).$
\end{defn}

In the following, let ${\cal F}_t$ denote the filtration generated by the random sequence $\{\phi_{[0,t]}\}$.
Define a sequence of stopping times $\{\stp_i : i \in \mathbb{N}_+\}$, measurable on the filtration described above, which is assumed to be non-decreasing, with $\stp_0=0$.

\begin{thm}
\label{thm5} \cite{YukMeynTAC2010}
Suppose that we have a $\varphi$-irreducible Markov chain $\phi$. Suppose moreover that there are functions $V \colon \mathbb{X} \to [a,\infty)$, $\beta \colon \mathbb{X} \to [a,\infty)$, $f\colon \mathbb{X} \to [a,\infty)$, for some $a \geq 0$, small set $C$, constant $b \in \Re$ and consider:
\begin{align}
&\Expect[V(\phi_{\stp_{z+1}}) \mid \clF_{\stp_z }]  \leq  V(\phi_{\stp_{z}}) -\beta(\phi_{\stp_z }) + b1_{\{\phi_{\stp_z} \in C\}}, \label{eq:thm5a}
\\
& \Expect \Bigl[\sum_{k=\stp_z}^{\stp_{z+1}-1} f(\phi_k)  \mid \clF_{\stp_z }\Bigr]  \le \beta(\phi_{\stp_z}), \quad  z\ge 0. \label{eq:thm5b}
\end{align}
If $a=1$ and \eqref{eq:thm5a} holds then $\phi$ is positive Harris recurrent with some unique invariant distribution $\pi$. If $a=0$, \eqref{eq:thm5a}, \eqref{eq:thm5b} hold and $\phi$ is positive Harris recurrent with some unique invariant distribution $\pi$ then we get that $ \lim_{t \to \infty} E[f(\phi_t)] < \infty. $
\end{thm}

\bibliographystyle{ieeetr}
\bibliography{SerdarBibliography}

\end{document}

%% file: symbols.tex
\def\mindex#1{\index{#1}}

\def\sq{\hbox{\rlap{$\sqcap$}$\sqcup$}}
\def\qed{\ifmmode\sq\else{\unskip\nobreak\hfil
\penalty50\hskip1em\null\nobreak\hfil\sq
\parfillskip=0pt\finalhyphendemerits=0\endgraf}\fi\medskip}

\long\def\defbox#1{\framebox[.9\hsize][c]{\parbox{.85\hsize}{%
\parindent=0pt
\baselineskip=12pt plus .1pt      
\parskip=6pt plus 1.5pt minus 1pt 
 #1}}}

\long\def\beginbox#1\endbox{\subsection*{}%
\hbox{\hspace{.05\hsize}\defbox{\medskip#1\bigskip}}%
\subsection*{}}

\def\endbox{}

\newsavebox{\junk}
\savebox{\junk}[1.6mm]{\hbox{$|\!|\!|$}}

\def\det{{\mathop{\rm det}}}
\def\limsup{\mathop{\rm lim\ sup}}

\newcommand{\field}[1]{\mathbb{#1}}

\def\Re{\field{R}}

\def\bfmath#1{{\mathchoice{\mbox{\boldmath$#1$}}%
{\mbox{\boldmath$#1$}}%
{\mbox{\boldmath$\scriptstyle#1$}}%
{\mbox{\boldmath$\scriptscriptstyle#1$}}}}

\def\bfmY{\bfmath{Y}}

\def\bfmhhaY{\bfmath{\hhaY}} 
\def\bfmhhaY{\hbox to 0pt{$\widehat{\bfmY}$\hss}\widehat{\phantom{\raise 1.25pt\hbox{$\bfmY$}}}}

\def\til={{\widetilde =}}

\def\clF{{\cal F}}

 \def\FRAC#1#2#3{\genfrac{}{}{}{#1}{#2}{#3}}

\def\ddtp{{\mathchoice{\FRAC{1}{d^{\hbox to 2pt{\rm\tiny +\hss}}}{dt}}%
{\FRAC{1}{d^{\hbox to 2pt{\rm\tiny +\hss}}}{dt}}%
{\FRAC{3}{d^{\hbox to 2pt{\rm\tiny +\hss}}}{dt}}%
{\FRAC{3}{d^{\hbox to 2pt{\rm\tiny +\hss}}}{dt}}}}

\def\Expect{E}

\def\average#1,#2,{{1\over #2} \sum_{#1}^{#2}}

\def\eye(#1){{\bf(#1)}\quad}

\newtheorem{theorem}{Theorem}[section]

\def\eq#1/{(\ref{e:#1})}

\newcommand{\beqn}[1]{\notes{#1}%
\begin{eqnarray} \elabel{#1}}

\newcommand{\eeqn}{\end{eqnarray} }

\newcommand{\beq}[1]{\notes{#1}%
\begin{equation}\elabel{#1}}

\newcommand{\eeq}{\end{equation}}

\def\bdes{\begin{description}}
\def\edes{\end{description}}

\newcounter{rmnum}

\newcounter{anum}

{\end{list}}

\def\ass(#1:#2){(#1\ref{#1:#2})}

\def\ritem#1{
\item[{\sf \ass(\current_model:#1)}]
}

\newenvironment{recall-ass}[1]{%
\begin{description}
\def\current_model{#1}}{
\end{description}
}

\newcommand{\bd}{\begin{description}}
\newcommand{\ed}{\end{description}}
\newcommand{\bt}{\begin{theorem}}
\newcommand{\et}{\end{theorem}}
\newcommand{\ba}{\begin{array}{rcl}}
\newcommand{\ea}{\end{array}}